\documentclass[leqno,11pt]{article}

\usepackage{amsfonts,amssymb,amsgen,theorem,stmaryrd, amsopn}
\usepackage[leqno]{amsmath} \usepackage[T1]{fontenc}
\usepackage[latin1]{inputenc} \usepackage{graphicx}
\DeclareGraphicsExtensions{.pdf,.ps}
\usepackage{graphpap}
\usepackage{times}

 \def\pasdegrille{\let\grille =
\pasgrille} \def\ecriture#1#2{\setbox1=\hbox{#1} \dimen1= \wd1
\dimen2=\ht1 \dimen3=\dp1 \grille #2 \box1 } \def\aat#1#2#3{ \divide
\dimen1 by 48 \dimen3=\dimen1 \multiply \dimen1 by #1 \advance \dimen1
by -\dimen3 \divide \dimen1 by 101 \multiply \dimen1 by 100 \divide
\dimen2 by \count11 \multiply \dimen2 by #2
\setbox0=\hbox{#3}\ht0=0pt\dp0=0pt \rlap{\kern\dimen1 \vbox
to0pt{\kern-\dimen2\box0\vss}}\dimen1= \wd1 \dimen2=\ht1}
\def\pasgrille{ \count12= \dimen1 \divide \count12 by 50 \divide
\dimen2 by \count12 \count11 =\dimen2 \ \divide \dimen1 by 48
\setlength{\unitlength}{\dimen1} \smash{\rlap{\ }} \dimen1= \wd1
\dimen2=\ht1 } \def\grille{ \count12= \dimen1 \divide \count12 by 50
\divide \dimen2 by \count12 \count11 =\dimen2 \ \divide \dimen1 by 48
\setlength{\unitlength}{\dimen1} \smash{\rlap{\graphpaper[1](0,0)(50,
\count11)}} \dimen1= \wd1 \dimen2=\ht1 }

\pasdegrille

\newtheorem{theoreme}{Theorem} \theorembodyfont{\sl}
\newtheorem{definition}{Definition}
\newtheorem{proposition}{Proposition} \newtheorem{lemme}{Lemma}
\newtheorem{rem}{Remark} \newtheorem{cor}{Corollary}
\theoremheaderfont{\sc}

 \newcommand{\Sg}{\mathrm{S}}

\renewcommand\Re{\mathrm{Re}\,} \renewcommand\Im{\mathrm{Im}\,}
\newcommand\R{{\mathbb R}}

\global
 \global
 \theoremstyle{break} \theorembodyfont{\it}
  \theoremheaderfont{\sc}

\newcommand{\dem}{\mbox{\it Proof: }} \newcommand{\cqfd}{\mbox{ }
\hfill$\Box$}

\newcommand{\BL}[4]{\dot{B}^{{#1},{#3}}_{#2}(\mathcal{L}^{#4}_t)}

\def\ra{\rightarrow} 
 \def\tq{\,\,\mbox{s.t.}\,\,}

\def\e{\varepsilon} \def\cdotv{\raise 2pt\hbox{,}}

\begin{document}

 \bibliographystyle{plain} \title{\vskip -1cm Smoothing And Dispersive
   Estimates For 1d Schr\"odinger Equations With Bv Coefficients And
   Applications}

  \author{Nicolas Burq \footnote{Département de mathématiques, UMR
  8628 du CNRS, B\^at 425 Université Paris-Sud, F-91405 Orsay} \ and Fabrice Planchon \footnote{ Laboratoire Analyse, G\'eom\'etrie
  \& Applications, UMR 7539, Institut Galil\'ee, Universit\'e Paris
  13, 99 avenue J.B. Cl\'ement, F-93430 Villetaneuse}}

 \date{} \maketitle
 \begin{abstract}
We prove smoothing estimates for Schr\"odinger equations \mbox{$i\partial_t \phi+\partial_x (a(x) \partial_x \phi)
=0$} with $a(x)\in \mathrm{BV}$, real and bounded from below. We then
bootstrap these estimates to obtain optimal Strichartz and maximal
function estimates, all of which turn out to be identical to the constant coefficient
case. We also provide counterexamples showing $a\in \mathrm{BV}$ to be
in a sense a minimal requirement. Finally, we provide an application to sharp
wellposedness for a generalized Benjamin-Ono equation.
\end{abstract}
 \par \noindent
\section*{Introduction}
Let us consider
\begin{equation}
  \label{eq:schrodbv}
  i\partial_t u+\partial_x(a(x)\partial_x u)=0,\,\,\,\,
  u(x,t=0)=u_0(x).
\end{equation}
We take $a\in \mathrm{BV}$, the space of bounded functions whose derivatives are Radon measures.
 Moreover, we assume $a$ to be real-valued and bounded from
below: $0<m\leq a(x)(\leq M)$.  We are interested in proving
smoothing and dispersive estimates for the function $u$. This type of equations has been
 recently studied by Banica~\cite{BaBV} who considered the case where
 the metric $a$ is piecewise constant (with a finite number of
 discontinuities).  In~\cite{BaBV}, Banica proved that the solutions
 of the Schr\"odinger equation associated to such a metric enjoy the
 same dispersion estimates (implying Strichartz) as in the case of the
 constant metric, and conjectured it would hold true for general $a\in
 \mathrm{BV}$ as well. Unfortunately, her method of proof (which
 consists in writing a complete description for the evolution problem)
 leads to constants depending upon the {\em number} of discontinuities
 rather than on the norm in $\mathrm{BV}$ of the metric and
 consequently does not extend to more general settings. On the other
 hand, Castro and Zuazua~\cite{CaZu} show that the space $\mathrm{BV}$
 is more or less optimal: they construct metrics $a\in C^{0,\beta}$
 for all $\beta \in [0,1[$ (but not in $\mathrm{BV}$) and solutions of
 the corresponding Schr\"odinger equation for which any {\em local}
 dispersive estimate of the type
$$\|u(t,x)\|_{L^1_{\text{loc}, t}( L^q_{\text{loc}, x})}\leq C
\|u_0\|_{H^s}$$ fail if $1/p < 1/2 -s$ (otherwise, the estimate is a
trivial consequence of Sobolev embeddings). In this article, we prove
the natural conjecture, namely that for $\mathrm{BV}$ metrics, the
Schr\"odinger equation enjoys the same smoothing, Strichartz and
maximal function estimates as for the constant coefficient case,
globally in time. In the context of variable coefficients, this
appears to be the first case where such a low regularity
(including discontinuous functions) is allowed, together with a translation
invariant formulation of the decay at infinity (no pointwise decay). Previous works on
dispersive estimates, while applying equally to any dimension, dealt
with $C^2$ compact pertubations of the Laplacian (\cite{StaTa}), short range pertubations with symbol-like decay at infinity
(\cite{RoZu}), and very recently long range pertubations, still with
symbol-like decay (\cite{HTW}). The idea to use local smoothing to derive Strichartz, however, goes back to
Staffilani-Tataru (\cite{StaTa}) in the context of variable
coefficients, and was used earlier to obtain full dispersion in \cite{JoSoSo} where a potential pertubation was
treated. All recent works on this topic make definitive use of resolvent
estimates for the elliptic operator, see
e.g. \cite{RodSch}. Finally, it has to be noticed that
Salort~\cite{Sa04} recently obtained dispersion (hence, Strichartz)
(locally in time) for 1D Schr\"odinger equations with $C^2$
coefficients through a completely different approach involving
commuting vector fields.

We now say a word on the relevance of non-trapping conditions.
 In higher dimension, it has been known since the works of Do\"
\i~\cite{Do96, Do00} and the first author~\cite{Bu02-2} that the non
trapping assumption is necessary for the optimal smoothing effect to
hold and the study of eigenfunctions on compact manifolds somewhat
shows that a non trapping condition is also necessary for Strichartz
estimates. In the one dimensional case, a smooth metric is always {\em
non trapping} as can be easily seen by a simple change of
variables. However, some trapping-related behaviours (namely the
existence of waves localized at a point) appear for metric with
regularity below $\mathrm{BV}$ (see the work by Castro
Zuazua~\cite{CaZu} and the appendix~\ref{sec.A2}). In fact the
assumption $a\in \mathrm{BV}$ ensures some kind of non trappingness
and this fact has been known for a while in the different context of
control theory~\cite{CoZu95}. Let us picture this on the model case of
piecewise constant metrics : consider a wave coming from minus
infinity. Then the wave propagates freely (at a constant speed) until
it reaches the first discontinuity. At this point some part of the
wave is reflected whereas some part is transmitted. It is easy to see
that a fixed amount of the energy (depending on the size of the jump
of the velocities) is transmitted. Then the transmitted wave propagate
freely until it reaches the second discontinuity, and so on and so
forth... Finally, we get that a fixed part of the energy of the
incoming wave is transmitted at the other end and propagates freely to
plus infinity. Whereas some part of the energy can remained trapped by
multiple reflections, this shows that some part is not trapped. As a
consequence, our geometry is (at least weakly) non trapping. This
phenomenon is clearly specific to the one dimensional case as can be
easily seen (using Snell law of refraction) with simple models
involving only two speeds.

The structure of our paper is as follows:
\begin{itemize}
\item In section \ref{sec1}, we prove a smoothing estimate which is
  the key to all subsequent results, by an elementary integration by
  parts argument, reminiscent of the time-space symmetry for the 1D
  wave equation. Transferring results from the wave to Schr\"odinger is
  sometimes called a transmutation and has been used in different
  contexts (\cite{Luc}).
\item We then obtain Strichartz and maximal function estimates in
  section \ref{sec2}, by combining our smoothing estimate with known
  estimates for the flat case.
\item Finally, we provide an application in section~\ref{BO},
  obtaining sharp wellposedness for a generalized Benjamin-Ono
  equation. Further applications to the Benjamin-Ono hierarchy of
  equations, including the true Benjamin-Ono, will be addressed
  elsewhere (\cite{BPBO}). The methods developed in this paper are likely to apply to other 1D dispersive models and quasilinear equations.
\item The first appendix is a short recollection of some results of
  Auscher-Tchamitchian~\cite{AT} and
  Auscher-MacIntosh-Tchamitchian~\cite{AMT} which imply that the
  spectral localization with respect to the operators $\partial_x a(x)
  \partial_x$ and $\partial_x ^2$ are reasonably equivalent.
\item In a second appendix we give a self-contained proof of a suitably modified version
of Christ-Kiselev Lemma (see~\cite{Christ-Kiselev}).
\item In a third appendix we prove that the BV regularity threshold
  is optimal in a different direction from \cite{CaZu}: there exist a
  metric $a(x)$ which is in $L^\infty \cap W^{s,1}$ for any $0\leq s
  <1$, bounded from below by $c>0$, and such that no smoothing effect
  nor (non trivial) Strichartz estimates are true (even with
  derivatives loss). This construction is very close in spirit to the
  one by Castro-Zuazua~\cite{CaZu}.
\end{itemize}

Besov spaces will be a convenient tool to state and prove many of our
results; we end this introduction by recalling their definition via
frequency localization (\cite{BL} for details).
\begin{definition}
\label{d1}
Let $\phi \in \mathcal{S}(\mathbb{R}^{n})$ such that $\widehat\phi =
1$ for $|\xi|\leq 1$ and $\widehat\phi= 0$ for $|\xi|>2$,
$\phi_{j}(x)= 2^{nj}\phi(2^{j}x)$, $ S_{j} = \phi_{j}\ast\cdot$,
$\Delta_{j} = S_{j+1} - S_{j}$. Let $f$ be in $\mathcal{S}'(\mathbb{
R}^{n})$.  We say $f$ belongs to $\dot B^{s,q}_{p}$ if and only if
\begin{itemize}
\item The partial sum $ \sum^{m}_{-m} \Delta_{j}(f)$ converges to $f$
as a tempered distribution (modulo polynomials if $s>n/p$ and $q>1$).
\item The sequence $\varepsilon_{j} = 2^{js}\| \Delta_{j} (f)\|_{L^{p}}$
belongs to $l^{q}$.
\end{itemize}  
\end{definition}
A suitable modification will be of interest, to handle the additional
time variable.
\begin{definition}
  \label{raah11}
Let $u(x,t)\in \mathcal{S}'(\mathbb{R}^{n+1})$, $\Delta_j$ be a
frequency localization with respect to the $x$ variable. We will say
that $u\in \BL s p q \rho$ iff
\begin{equation}
  \label{eq:raaah12}
  2^{js}\|\Delta_j u\|_{L^p_x(L^{\rho}_t)} =\varepsilon_j \in l^q,
\end{equation}
and other requirements are the same as in the previous definition.
\end{definition}
Notice that whenever $q=\rho$, the Besov space $\BL s p q q$ is
nothing but the usual ``Banach valued'' Besov space $\dot
B^{s,q}_p(F)$ with $F=L^q_t$.

 Finally, through this article we will denote by $\Sg (t) = e^{it
 \partial_x ^2}$ and $\Sg_a(t) = e^{it \partial _x a(x) \partial_x}$
 the (1D) group-evolution associated to the constant and variable coefficients equations respectively.

{\bf Acknowledgments:} we would like to thank Pascal Auscher, Isabelle Gallagher and Christian G\'erard for discussions on various aspects of this paper.

\section{Local smoothing}\label{sec1}
For the (flat) Schr\"odinger equation on the real line, we have the
following estimate:
$$ \|\partial_x \Sg (t) \phi_0\|_{L^\infty_x L^2_t}\simeq
\|\phi_0\|_{\dot H^{\frac 1 2}}.
$$ It can be proved directly using the Fourier transform (see
\cite{KPV91}). With this in mind, one can of course write a similar
estimate for the 1D wave equation, which is also a
trivial consequence of the explicit representation as a sum of
traveling waves; however, one can prove it as well by integration by parts on
the inhomogeneous equation, exchanging $t$ and $x$ which play
equivalent roles. This last procedure is flexible enough to allow
variable coefficients and will lead to our first result. We start by
stating once and for all our hypothesis on the coefficient $a$.
\begin{definition}
  We call $a$ an m-admissible coefficient when the following
  requirements are met:
  \begin{itemize}
  \item the function $a$ is real-valued, belongs to $\mathrm{BV}$,
    namely
$$ \partial_x a \in \mathcal M=\{ \mu \tq \int_\R d|\mu|<+\infty\},
$$
\item the function $a$ is bounded from below almost everywhere by $m$.
  \end{itemize}
 We will denote by $M$ its maximum and
  $\|a\|_{\mathrm{BV}}$ its bounded variation ($a(x)\leq M\leq \|a\|_{\mathrm{BV}}$).
\end{definition}

After this preliminary definition, we can state the main theorem.
\begin{theoreme}
\label{th:smoothing}
Let $m>0$ and $a$ be an m-admissible coefficient. There exist
  $C(\|a\|_{\mathrm{BV}}, m)>0$ such that
  \begin{itemize}
  \item If $u,f$ are solutions of
\begin{equation}\label{eq.nonhomog}
(i\partial_t +\partial _x a(x) \partial _x) u =f,
\end{equation}
with zero Cauchy data then
\begin{equation}\label{eq:visco}
\|\partial _x u \|_{L^\infty_x L^2_t}+ \| (-\partial^2_t)^{1/4} u
\|_{L^\infty_x L^2_t} \leq C \|f\|_{L^1_x L^2_t}.
 \end{equation}
\item If
$$(i\partial_t +\partial _x a(x) \partial _x ) u=0,\text{ with }
u_{|t=0} = u_0 \in L^2$$ then
\begin{equation}
  \label{eq:smoothing-homogene}
\| u\|_{\BL {\frac 1 2} \infty 2 2}\leq C \|u_0\|_{L^2}.
\end{equation}
  \end{itemize}
\end{theoreme}
\begin{rem}
  One may wonder why we chose to consider $\partial_x a(x)\partial_x$
  as opposed to, say, $g(x)\partial^2_x$. It turns out that one may
  obtain one from another through an easy change of variable, and we
  elected to keep the divergence form as the most convenient for
  integration by parts. The astute reader will check that
  $b(x)\partial_x a(x)\partial_x$ can be dealt with as well, and the
  additional requirement will be for $b$ to be $m$-admissible. Remark
  also that our method can handle first order terms of the kind
  $b(x)\partial_x$ with $b\in L^1$ (see section \ref{BO}).

\end{rem}
\dem In order to obtain (\ref{eq:visco}), we will reduce ourselves to
a situation akin to a wave equation and perform an integration by
parts. Obtaining (\ref{eq:smoothing-homogene}) from (\ref{eq:visco})
is then a simple interpolation and $TT^\star$ argument. We first reduce the study to smooth $a$.
\begin{proposition}
\label{regul} Denote by $A= \partial_x a(x) \partial_x$. Assume that the evolution
semi-group $\Sg_a(t)$ satisfies for any smooth ($C^\infty$)
m-admissible $a$:
$$
 \forall u_0 \in L^2, \qquad \|\Sg_a(t) u_0\|_B\leq C \|u_0\|_{L^2}
$$
 with $B$ a Banach space (weakly) continuously embedded in
$\mathcal{D}'(\mathbb{R}^2)$, whose unit ball is weakly compact, and
$C$ a constant depending only on $m$ and $\|\partial_x
a\|_{L^1}$. Then the same result holds (with the same constant) for
any m-admissible $a$.
\end{proposition}
\dem Let us consider $\rho\in C^\infty_0( \mathbb{R})$ a non negative
function such that $\int \rho=1$, and $\rho_\varepsilon=
\varepsilon^{-1}\rho(x/\varepsilon)$. Denote by $a_\varepsilon=
\rho_\varepsilon \star a$ and $A_\e=-\partial_x
a_\e(x)\partial_x$. The sequence $a_\varepsilon$ is bounded in
$\dot{W}^{1,1}$. Furthermore, $a_\varepsilon$ converges to $a$ for the
$L^\infty$ weak $\star$ topology.  According to the weak compactness
of the unit ball of $B$, taking a subsequence, we can assume that
$\Sg_{a_\varepsilon}(t) u_0$ converges weakly to a limit $v$ in $B$ (and
consequently in $\mathcal{D}'( \mathbb{R}^2)$). To conclude, it is
enough to show that $v= \Sg_a(t) u_0$ in $\mathcal{D}'(
\mathbb{R}^2)$.  We first remark that  as a (multiplication) operator on
$L^2$, $a_\varepsilon$ converges strongly to $a$ (but of course not in
operator norm) and consequently $\partial_x a_\varepsilon(x) \partial_x$
converges strongly to $\partial_x a(x) \partial_x$ as operators from
$H^1$ to $H^{-1}$. On the other hand the bound $0<m\leq a(x) \leq M$
and the fact that $\rho$ is non negative imply that $a_\varepsilon$
satisfy the same bound and consequently that the family $(A_\varepsilon +
i)^{-1}$ is bounded from $H^{-1} $ to $H^1$ by $1/m$.  From the
resolvent formula
$$ (A_\varepsilon + i)^{-1} - (A+i)^{-1} = (A_\varepsilon + i)^{-1}(
A - A_\varepsilon)(A + i)^{-1},$$ given $(A_\varepsilon +i)^{-1}$ is
uniformly bounded from $H^{-1}$ to $H^1$, we obtain that $ (A_\varepsilon
+ i)^{-1}$ converges strongly to $(A+i)^{-1}$ as an operator from
$H^{-1}$ to $H^1$, and consequently as an operator on $L^2$. This
convergence implies (see~\cite[Vol I, Theorem VIII.9]{ReSi78}) that
$A_\varepsilon$ converges to $A$ in the strong resolvent sense and
(see~\cite[Vol I, Theorem VIII.21]{ReSi78}) that for any $t\in
\mathbb{R}$, $\Sg_{a_\varepsilon}(t)$ converges strongly to
$\Sg_a(t)$. Finally, from the boundedness of $\Sg_{a_\varepsilon}(t)u_0$
in $L^\infty_t( L^2_x)$, we deduce by dominated convergence that
$\Sg_{a_\varepsilon}(t) u_0$ converges to $\Sg_{a}(t) u_0$ in
$L^1_{t,\text{loc}}(L^2)$ and hence in $\mathcal{D}'$. Similarly, we
can handle non-homogeneous estimates.\cqfd

\begin{rem}
  Alternatively, we can perform our argument for $a$ a step function
  with finite $\mathrm{BV}$ norm. We will briefly sketch this at the
  end of this section.
\end{rem}
We are now considering the following equation (for $a\in C^\infty_0$):
\begin{equation}
  \label{eq:resolvent}
 -\sigma v+\partial_x (a(x)\partial_x v)=g.
\end{equation}
where $v, g$ will be chosen later to be the time Fourier transform of
$u, f$. 
\begin{proposition}\label{prop.res}
There exist $C(m, \|a\|_{\mathrm{BV}})$ such that for any $\sigma =
\tau + i \varepsilon, \varepsilon \neq 0$ the resolvent $ (-\sigma
+\partial_x a(x)\partial_x)^{-1}$, which is a well defined operator
from  $L^1\subset H^{-1}$ to $H^1\subset L^\infty$ and from $L^2$ to $H^2$ satisfies
\begin{equation}
  \label{eq:smooth-resolvent}
 \|(-\sigma +\partial_x a(x)\partial_x)^{-1}\|_{L^1 \rightarrow
 L^\infty}\leq C .
\end{equation}
\end{proposition}
It should be
noticed that since this and all further estimates are scale invariant (including the
constants which are dependent on scale invariant quantities of $a$),
we could reduce the study to the case $\tau=\pm 1$ by changing $a(x)$
into $a(\sqrt {\pm\tau}^{-1} x)$. We elected to keep $\tau$ through
the argument as it helps doing book keeping.
\begin{rem}
The elliptic case ($\tau>0$) is more or less understood and as a
corollary, the associated heat equation as well. In fact these results
apply to a larger class of $a$ than the one we consider here: $a\in
L^\infty$, $\Re a>0$. More specifically, the heat kernel (and its
derivatives) associated to the operator
$A=-\partial_x(a(x)\partial_x)$ is known to be of Gaussian type, a
fact which will be of help to handle derivatives. A very nice and
thorough presentation of this (and a lot more !) can be found in
\cite{AMT}. We refer to Appendix \ref{equiv-norm} for a short
recollection of the facts we will need later.
\end{rem} 
In the sequel we will perform integrations by parts. We can assume $g\in L^2$. Consequently $v\in H^2$ and these integrations by parts are licit (in
particular, the boundary terms near $\pm \infty$ vanish).  We first
multiply \eqref{eq:resolvent} by $\overline{v}$, integrate by parts
and take the imaginary and real parts.  This yields
\begin{equation}\label{eq.apriori1}
\begin{gathered}
|\varepsilon| \int_{\R} |v| ^2 \leq \|g \|_{L^1}\|v\|_{L^\infty}\\
|\varepsilon| \int_{\R} a(x)|\partial_x v| ^2 \leq |\varepsilon| |\tau|
\int_{\R} |v| ^2 + |\varepsilon| \|g\|_{L^1}\|v\|_{L^\infty}\leq
(|\varepsilon|+ |\tau|) \|g \|_{L^1}\|v\|_{L^\infty}
\end{gathered}
\end{equation}
We now proceed in the hyperbolic region $-\tau>0$.  Multiplying
(\ref{eq:resolvent}) by $a(x)\partial_x \overline{v}$ and integrating,
we get
\begin{equation}
\label{eq:hyper-ell}
\int_{-\infty}^x -\sigma a v \partial_x(\overline{ v})+
\int_{-\infty}^x \partial_x (a\partial_x v) a\partial_x \overline{v}
= \int_{-\infty}^{x} g a\partial_x \overline{v}.
\end{equation}
Integration by parts and taking the real part yields
\begin{multline}
  \label{eq:hyper}
  -\tau a|v|^2(x)+ |a\partial_x v|^2(x) +2 \int_{-\infty}^x \tau
(\partial_x a)|v|^2 \\ \leq 2|\varepsilon| \int_{\R} a |v| |\partial_x
v| +2\|g\|_{L^1} \|a\partial_x v\|_{L^\infty} .
\end{multline}
We now use~\eqref{eq.apriori1} to estimate the right hand side
in~\eqref{eq:hyper} and obtain
\begin{multline*}
  -\tau a|v|^2(x)+ |a\partial_x v|^2(x) + 2\int_{-\infty}^x \tau
(\partial_x a)|v|^2 \\ \leq 2\max(1,
\|a\|^{\frac 1 2}_{L^\infty})\|g\|_{L^1} \bigl( \|a\partial_x v\|_{L^\infty}+
(|\varepsilon|+|\tau|)^{1/2}\|v\|_{L^\infty}\bigr).
\end{multline*}
On the other hand we are in 1D and,
\begin{equation}\label{eq.1D}
 \|v\|_{L^\infty}^2 \leq 2 \|v\|_{L^2} \|\partial_xv\|_{L^2}
\end{equation}
which implies, using~\eqref{eq.apriori1},
\begin{equation*}
\varepsilon\|v\|_{L^\infty}^2 \leq 2 m^{-\frac 1 2} \|g\|_{L^1} \sqrt{ |\varepsilon| +
|\tau|} \|v\|_{L^\infty}
\end{equation*}
Consequently we get
\begin{multline*}
  (|\varepsilon|+ |\tau|) a|v|^2(x)+ |a\partial_x v|^2(x) +
2\int_{-\infty}^x \tau (\partial_x a)|v|^2 \\ \leq
C(m,\|a\|_{\mathrm{BV}})\|g\|_{L^1} \bigl( \|a\partial_x v\|_{L^\infty}+
(|\varepsilon|+|\tau|)^{1/2}\|v\|_{L^\infty}\bigr).
\end{multline*}
Setting
\begin{eqnarray*}
\Omega_+(x) & = &\sup_{y< x} \left(|\varepsilon|+|\tau|) a(y)
 |v|^2(y)+|a(y)\partial_x v|^2(y)\right)\\ k(x) & = & a(x)^{-1}
 |\partial_x a|,
\end{eqnarray*}
we have
\begin{equation}\label{eq.gron}
 \Omega_+(x)\leq C(m, \|a\|_{\mathrm{BV}}) \sqrt \Omega_+(+\infty)
\|g\|_{L^1_x}+2\int_{-\infty}^x k (y) \Omega_+(y) dy.
\end{equation}
 Given that $\Omega_+$ is positive, we obtain by Gronwall inequality
\begin{equation*}
\begin{aligned}
\int_{-\infty} ^x k(y) \Omega_+(y) dy &\leq C(m, \|a\|_{\mathrm{BV}})
\left(\int_{-\infty} ^x e^{ \int_{y} ^x 2k(z)dz } k(y)
dy\right)\|g\|_{L^1_x} \sqrt \Omega_+(+\infty)\\ &\leq 2 C(m,
\|a\|_{\mathrm{BV}}) e^{\int_{-\infty}^x 2k(y)dy}\|g\|_{L^1_x}\sqrt
\Omega_+(+\infty)
\end{aligned}\end{equation*}
and consequently, coming back to~\eqref{eq.gron}
\begin{equation}
 \sqrt{\Omega_+(+\infty)}\leq C (m, \|a\|_{\mathrm{BV}})\|g\|_{L^1}\bigl( 2+
 8e^{2\|k(x)\|_{L^1}}\bigr).
\end{equation}
Now we proceed with the elliptic region $\tau>0$, for which the above
line of reasoning fails.  We perform the usual elliptic regularity
estimate and multiply the equation by $\overline{v}$, to obtain
\begin{equation*}\label{eq:elli}
 \int_{\R}\tau |v|^2+a|\partial_x v|^2 =-\Re \int_{\R} g
\overline{v}\, , \qquad \varepsilon \int_{\R} |v|^2 = - \Im \int_{\R} g
\overline{v}
\end{equation*}
 which gives
\begin{equation}\label{eq.estmell2}
\int_{\R} \left((|\tau| + |\varepsilon|) |v|^2+a|\partial_x v|^2\right)
\leq 2\| g\|_{L^1_x} \|v\|_{L^\infty } .
\end{equation}
In order to conclude, we go back to the (beginning of) the estimate we
made in the hyperbolic case, i.e. (\ref{eq:hyper-ell}) and integrate
by parts only the second term in the left hand side,
\begin{equation*}
|a\partial_x v|^2(x)\leq 2\int_{-\infty}^x |g| a|\partial_x v| {}+ 2
\int_{-\infty}^x |\sigma| a |v||\partial_x v|
\end{equation*}
and to bound the last term we use~\eqref{eq.estmell2},
\begin{equation}\label{eq.estimell}
\|a\partial_x v\|_{L^\infty}^2 \leq \|g\|_{L^1_x}\bigl(2 \|a
\partial_x v\|_{L^\infty} + 4 |\tau| ^{1/2}\|v\|_{L^\infty}\bigr).
\end{equation}
Adding $\tau a|v|^2$ to~\eqref{eq.estimell} and
using~\eqref{eq.estmell2},~\eqref{eq.1D}, we obtain
\begin{equation*}
  \begin{aligned}
\Omega_-(x) & = \sup_{y\leq x} (|\varepsilon|+|\tau|) a|v|^2
(y)+|a\partial_y v|^2(y)\\
 & \leq 2M |\tau \|v\|_{L^2}
\|\partial_xv\|_{L^2} +4
(|\varepsilon|+|\tau|)^{1/2}\|g\|_{L^1_x}\|v\|_{L^\infty}\\
 & \leq
\|g\|_{L^1} (2M +4)\|(|\varepsilon|+|\tau|)^{1/2}v\|_{L^\infty}
  \end{aligned}
\end{equation*}
which gives again
\begin{equation}\label{eq:resestbis}
\sup_x \Omega_-(x)\leq \frac{(\|a\|_{L^\infty}+4)^2} m \|g \|^2_{L^1}.
\end{equation}
This ends the proof of Proposition~\ref{prop.res}.\cqfd
\begin{rem}
  Notice that for this elliptic estimate, we only used $a\in L^\infty$
  and nothing else.
\end{rem}
We now come back to the proof of Theorem~\ref{th:smoothing}. Consider
$u,f$ solutions of~\eqref{eq.nonhomog}. We can assume that $f$ (and
consequently $u$) is supported in $t>0$ (the contribution of negative
$t$ being treated similarly)). Then for any $\varepsilon>0$ $u_\varepsilon
= e^{-\varepsilon t} u$ is solution of
$$ (i\partial_t +i\varepsilon + \partial_x a(x) \partial_x) u_\varepsilon
=f,
$$ Assuming that $f$ has compact support (in time), we can consider
the Fourier transforms with respect to $t$ of $f$ and $u_\varepsilon$, $g(\tau)$ and
$v_\varepsilon(\tau)$ which satisfy
$$ (- \tau +i\varepsilon + \partial_x a(x) \partial_x) v_\varepsilon =g.
$$ We may now apply Proposition~\ref{prop.res}, take $L^2_\tau$ norms,
switch norms and revert back to time by Plancherel, and get
\begin{equation*}
\begin{aligned}
 \|\partial_x
u_\varepsilon\|_{L^\infty_x(L^2_t)}+\|(-\partial^2_t)^{\frac 1 4}
u_\varepsilon\|_{L^\infty_x(L^2_t)} & =\|\partial_x
v_\varepsilon\|_{L^\infty_x(L^2_\tau)}+\|(-\partial^2_t)^{\frac 1 4}
v_\varepsilon\|_{L^\infty_x(L^2_\tau)}\\
 & \leq \|\partial_x
v_\varepsilon\|_{L^2_\tau (L^\infty_x)}+\|(-\partial^2_t)^{\frac 1 4}
v_\varepsilon\|_{L^2_\tau (L^\infty_x)}\\
&  \leq C \|g_\varepsilon\|_{L^2_\tau
(L^1_x)} \leq C\|g_\varepsilon\|_{L^1_x(L^2_\tau)}=
C\|f_\varepsilon\|_{L^1_x(L^2_t)},
\end{aligned}
\end{equation*}
where $C= C(m,\|\partial_x a\|_{L^1_x})$ is uniform with respect to
 $\varepsilon >0$. Letting $\varepsilon>0$ tend to $0$, we obtain the same
 estimate for $u$, which is exactly \eqref{eq:visco} in Theorem
 \ref{th:smoothing} (up to replacement of $\mathrm{BV}$ by $\dot
 W^{1,1}$, which was dealt with in Proposition \ref{regul}). Finally
 we easily drop the compact in time assumption for $f$ by a density argument.

We are left with proving the homogeneous estimate
\eqref{eq:smoothing-homogene}. As usual, estimates on the homogeneous
problem follow from the estimate with a fractional time derivative: by
a $TT^\star$ argument, and using the commutation between time
derivatives and the flow, we get
$$ \|(-\partial^2_t)^{\frac 1 8} u\|_{L^\infty_x(L^2_t)}\lesssim \sqrt
C \|u_0\|_{L^2_x}.
$$ Then, using the equation, $i\partial_t u=Au$ where $A=-\partial_x
a(x) \partial_x$, we can replace $(i\partial_t)^{1/4}$ by
$A^{1/4}$. However, we will need real derivatives later, rather than
powers of $A$. We postpone the issue of equivalence between the two
and take another road: notice that we obtained
\eqref{eq:smooth-resolvent} for solutions of \eqref{eq:resolvent}
$$ \|\partial_x v\|_{L^\infty_x}\lesssim \|g\|_{L^1_x},
$$ which immediately implies 
 \begin{equation}
   \label{eq:besov1}
\| v\|_{\dot B^{1,\infty}_\infty}\lesssim \|g\|_{\dot B^{0,1}_1}.
 \end{equation}
Call $R_\sigma= (\partial_x a(x) \partial_x - \sigma ) ^{-1}$. Its
adjoint is $R_{\overline \sigma}$ and according to
Proposition~\ref{prop.res} (applied to $\overline \sigma = \tau - i
\varepsilon$), we get
\begin{equation}
  \label{eq:besov-1}
 \| v\|_{\dot B^{0,\infty}_\infty}\lesssim \|g\|_{\dot B^{-1,1}_1}.
\end{equation}
 By real interpolation (recall $\left(\dot B^{s_1,q_1}_p,\dot
B^{s_2,q_2}_p\right)_{\theta,r}=\dot B^{s,r}_p$), we obtain (with
$\theta=1/2, r=2$)
$$ \| v\|_{\dot B^{\frac 1 2,2}_\infty}\lesssim \|g\|_{\dot B^{-\frac
1 2,2}_1}.
$$ Given that the third index is $2$, we can again take $L^2_\tau$
 norms, switch them (Minkowski) and by Plancherel (and letting $\varepsilon$ tend to $0$), we get the desired
 estimate:
\begin{equation*}
  \label{eq:besov-inhomogene}
\| u\|_{\dot B^{\frac 1 2,2}_\infty(L^2_t)}\lesssim \|f\|_{\dot
B^{-\frac 1 2,2}_1(L^2_t)}.
\end{equation*}
 denote by $\Sg_a(t)$ the evolution group for the homogeneous
equation, we have
$$ u=\int_{s<t} \Sg_a(t-s) f(s)ds,
$$ solution of the inhomogeneous problem, and we can as well treat the
$s>t$ case. Hence we have obtained
$$ \|\int \Sg_a(t-s)f(s)ds\|_{\dot B^{\frac 1 2,2}_2( L^2_t)}\lesssim
 \|f\|_{\dot B^{-\frac 1 2,2}_1 (L^2_t)}.
$$ The usual $TT^\star$ argument applies and gives
$$ \|\Sg_a(t) u_0\|_{\dot B^{\frac 1 2,2}_\infty (L^2_t)}\lesssim
 \|u_0\|_{L^2_x}.
$$ This ends the proof of Theorem \ref{th:smoothing}.\cqfd

We now provide an alternative argument which directly proves the
resolvent estimate for $a$ a step function, bounded from below and
with bounded variation. For the sake of conciseness, we take directly
$\sigma= \tau \in \mathbb{R}$ and will not justify the validity of the
integration by parts (and in particular the vanishing of the boundary
terms at $\pm \infty$). As before, the justification consists in taking
$\sigma = \tau + i \varepsilon$ and passing to the limit $\varepsilon
\rightarrow 0$. We set, $m=1$ for simplicity, and rescale to obtain
$\tau=\pm 1$. Starting from \eqref{eq:resolvent}, with $\tau=+1$ (the
difficult case) and denoting
$$ a(x)=\sum_i a_i \chi_{[x_i,x_{i+1}[}(x), \text{ with } \sum_i
|a_i-a_{i-1}|<+\infty,
$$ we have, with $x\in [x_I,x_{I+1}[$
$$ \sum_{i<I} \int_{x_i}^{x_{i+1}} a_i v\partial_x \bar v+\int_{x_I}^x
a_I v\partial_x \bar v+a_I^2 |\partial_x v|^2(x)=\int_{-\infty }^x g
a\partial_x \bar v.
$$ Integrating by parts,
$$ \sum_{i<I} a_i (|v|^2(x_{i+1})-|v|^2(x_{i}))+a_I
(|v|^2(x)-|v|^2(x_I))+a_I^2 |\partial_x v|^2(x)\leq 2 \| g\|_1
\sup_{y\leq x} a(y)| \partial_x v(y)|.
$$
$$ \sum_{i-1<I} (a_{i-1}-a_i) |v|^2(x_{i})+ |v|^2(x)+a^2(x)
|\partial_x v|^2(x)\leq 2 \| g\|_1 \sup_{y\leq x} a(y) | \partial_x
v(y)|.
$$ We rewrite this as
\begin{equation}
  \label{eq:contdis}
\sup_{y\leq x}( 2 |v|^2(y)+|\partial_x v|^2(y))\leq 4 \|
g\|^2_1+\sum_{i\leq I} |a_{i-1}-a_i| |v|^2(x_{i})
\end{equation}
at which point one may simply replace the left hand side (noting that
$x<x_{I+1}$) by its weaker discrete counterpart
$$ \sup_{i\leq I+1}( 2 |v|^2(x_i)+|\partial_x v|^2(x_i))\leq 4 \|
g\|^2_1+\sum_{i\leq I} |a_{i-1}-a_i| |v|^2(x_{i})
$$ which becomes a discrete analog of Gronwall: call
$\gamma_i=\sup_{j\leq i} ( 2 |v|^2(x_j)+|\partial_x v|^2(x_j))$ and
$\alpha_i=|a_{i-1}-a_{i}|$, we have
$$ \gamma_{I+1}\leq C+\sum_{i\leq I} \alpha_i \gamma_i.
$$ Therefore,
$$ \sum_{i\leq I} \frac{\alpha_{i+1} \gamma_{i+1}}{C+\sum_{j\leq i}
\alpha_j \gamma_j} \leq \sum_{i\leq I+1}\alpha_i<+\infty,
$$ Call $S_I=\sum_{i\leq I} \alpha_i \gamma_i$,
$$ \sum_{i\leq I}\int_{S_i}^{S_{i+1}} \frac {dx}{1+x} \leq \sum_{i\leq
I} \frac{ S_{i+1}-S_{i}}{1+S_i}\leq \sum_i \alpha_i,
$$ yielding
$$ S_I \lesssim \exp(\sum_i \alpha_i),
$$ which is nothing but the desired bound: recall \eqref{eq:contdis}
and notice we just bounded the right hand side.\cqfd

Notice that up to this point we avoided to use any of the machinery
presented in Appendix \ref{equiv-norm}, thus keeping the proof
self-contained. However, a rather natural question is now how one can
handle (fractional) derivatives: i.e., replace $u_0\in L^2$ by
$u_0\in \dot H^s$. In order to deal with commutation, we will rely in
a very natural way on Appendix \ref{equiv-norm}.
\begin{proposition}
\label{th:smoothing-s}
Assuming $a$ is m-admissible, we have:
\begin{itemize}
\item if $u,f$ are solutions of
\begin{equation*}
(i\partial_t +\partial _x a(x) \partial _x) u =f,
\end{equation*}
then, for $0< s<1$,
\begin{equation}
\label{eq:visco-s}
\| u \|_{\BL s \infty 2 2}\lesssim\|f\|_{\BL {s-1} 1 2 2}.
 \end{equation}
\item if $-1 < s <\frac 1 2$ and
$$ (i\partial_t +\partial _x a(x) \partial _x) u=0,\text{ with }
u_{|t=0} = u_0
$$ then
\begin{equation}
  \label{eq:smoothing-homogene-s}
\| u\|_{\BL {s+\frac 1 2} \infty 2 2 }\lesssim \|u_0\|_{\dot H^s}.
\end{equation}
\end{itemize}
\end{proposition}
\dem Recall that by interpolation between \eqref{eq:besov1} and
\eqref{eq:besov-1} we have
$$ \| v\|_{\dot B^{s,2}_\infty} \lesssim \| g\|_{\dot B^{s-1,2}_1},
$$ for all $0< s< 1$, which immediately gives \eqref{eq:visco-s}. For
the homogeneous problem, we simply rely on the equivalence
properties stated in Appendix \ref{sec:besov-eq}: we apply
\eqref{eq:smoothing-homogene} to $\Delta^A_j u_0$, a datum localized
with respect to $A$ (see the Appendix for a definition) and use
commutation between $\Delta_j^A$ and $S_a(t)$ to obtain
$$ \|\Delta^A_j u\|_{\BL {\frac 1 2} \infty 2 2} \lesssim \|\Delta^A_j
u_0 \|_{L^2_x}.
$$
Equivalence between $\BL {\frac 1 2} \infty 2 2$ and $\dot
B^{\frac 1 2,2}_{\infty,A} (\mathcal{L}^2_t)$ yields
$$ 2^{\frac 1 2 j}\|\Delta^A_j u\|_{L^\infty_x L^2_t} \lesssim
\|\Delta^A_j u_0 \|_{L^2_x},
$$ for which multiplying by $2^{js}$ and summing over $j$ provides the
desired result, after switching back from $A$ based Besov spaces to
the usual ones. Hence $s>-1$ from the right hand side, and $s+1/2<1$
from the left hand side.\cqfd

\section{Strichartz and maximal function estimates}\label{sec2}
We now prove Strichartz and maximal function estimates by
combining the smoothing effect from the previous section with a change
of variable and corresponding estimates for the flat Schr\"odinger equation.
\begin{theoreme}
\label{Strichartz-L2}
  Let $a$ be an m-admissible coefficient. Let $u$ be a solution of
  (\ref{eq:schrodbv}) with $u_0\in L^2$. Then for $\frac{2}{p}+\frac 1
  q=\frac 1 2$, $p\geq 4$, we have
  \begin{equation}
    \label{eq:Strichartz}
\|\Sg_a(t) u_0\|_{L^p_t(\dot B^{0,2}_q)} \lesssim \|u_0\|_{L^2}.
  \end{equation}
When $p>4$ ($q<+\infty$),
\begin{equation}
  \label{eq:Strichartz-nep}
      \|\Sg_a(t) u_0\|_{L^p_t(L^q_x)} \lesssim \|\Sg_a(t)
      u_0\|_{L^p_t(\dot B^{0,2}_q)} \lesssim \|u_0\|_{L^2}.
\end{equation}
\end{theoreme}
\begin{rem}
Notice that the end-point $(4,\infty)$ is missing. This can be seen as
an artifact of the proof. It will be clear that in this section, we
only use $a\in L^\infty\cap \dot B^{1,\infty}_1$ and bounded from below
(together with the estimates of Theorem \ref{th:smoothing}). Adding a
technical hypothesis like $a\in \dot B^{1,2}_1$ (which does not follow
from $a\in \mathrm{BV}$) would allow to recover the end-point, at the
expense of extra technicalities which we elected to keep out (see
\cite{BPBO} for further developments).
\end{rem}
One may state a corollary including fractional derivatives as well.
\begin{proposition}
\label{Strichartz-Hs}
  Let $u$ be a solution of (\ref{eq:schrodbv}), and $u_0\in \dot H^s$,
  $|s|< 1$. Then for $\frac{2}{p}+\frac 1 q=\frac 1 2$, $p\geq 4$, we
  have
  \begin{equation}
    \label{eq:Strichartz-frac}
    \|\Sg_a(t) u_0\|_{L^p_t(\dot B^{s,2}_q)} \lesssim \|u_0\|_{\dot
    H^s}.
  \end{equation}
\end{proposition}

Similarly, we also obtain maximal function estimates.
\begin{theoreme}
\label{maximal}
  Let $u$ be a solution of (\ref{eq:schrodbv}), and $u_0\in \dot H^s$,
  $-3/4< s<1$. Then
  \begin{equation}
    \label{eq:maximal}
    \|\Sg_a(t) u_0\|_{\BL{s-\frac 1 4} 4 2 \infty} \lesssim
    \|u_0\|_{\dot H^s}.
  \end{equation}
\end{theoreme}
\dem We aim at taking advantage of an appropriate new formulation for
our original problem and proving Theorems \ref{Strichartz-L2} and
\ref{maximal} at once. The operator $\partial_x a \partial_x$ may be
rewritten as $(\sqrt a \partial_x)^2+(\partial_x \sqrt a) \partial_x$,
and one would like to ``flatter out'' the higher order term through a
change of variable. However, performing directly a change of variable
leads to problems when dealing with the newly appeared first order
term. Therefore, we have to paralinearize the equation. Let us rewrite
$a$:
$$ a=\frac m 2+b^2,\text{ with } \partial_x b \in L^1_x,
$$ given that $a$ is $m$-admissible.  Writing
\begin{eqnarray*}
b^2 \partial_x u & = & \sum_k (S_{k-3} b)^2 \partial_x S_k u-(S_{k-4} b)^2
\partial_x S_{k-1} u\\ & = & \sum_k (S_{k-3} b)^2 \partial_x \Delta_k u +\sum_k \Delta_{k-3}
b (S_{k-3}+S_{k-4}) b \partial_x S_{k-1} u
\end{eqnarray*}
and applying $\Delta_j$ to the equation,
$$ i\partial_t \Delta_j u+\frac m 2 \partial^2_x\Delta_j u+\Delta_j
 \partial_x \sum_{k\sim j} \left((S_{k-3} b)^2\partial_x \Delta_k
 u\right)+\Delta_j \partial_x \sum_{j\lesssim k\sim l} \left(\Delta_{k-3}
 b S_{l-3} b \partial_x S_{k-1} u\right)= 0.
$$ From now on, we ignore shifts in indices for the last term as they
won't play any role.  Thus we get
$$ i\partial_t \Delta_j u+\left(\sqrt{\frac m 2+(S_{j-3} b)^2}\right)
 \partial_x \left( \left( \sqrt{\frac m 2+(S_{j-3} b)^2}\right)
 \partial_x \Delta_j u\right)=R_j,
$$ and, with $\tilde\Delta_j$ an enlargement of the localization,
 \begin{multline*}
    R_j=-\Delta_j \partial_x \sum_{j\lesssim k} \Delta_k b S_k b
\partial_x S_k u-\tilde\Delta_j \partial_x \sum_{k\sim j}
[\Delta_j,(S_{k-3}b)^2] \partial_x \Delta_k u\\ {}-S_{j-3} b
(\partial_x S_{j-3} b) ( \partial_x \Delta_j u).
 \end{multline*}
Assuming the smoothing effect from Theorem \ref{th:smoothing}, we can
effectively estimate the reminder.
\begin{proposition}\label{reminder}
  Assume that the hypothesis of Theorem \ref{th:smoothing} hold: then
 $$\ \sum_j R_j \in \BL {-\frac 1 2} 1 2 2.
$$
\end{proposition}
\dem Let us do the first term:  relabeling $k=j$ for simplicity,
$$ 2^{-\frac j 2} \partial_x S_j u \in l^2_j L^\infty_x L^2_t,
\,\,\,S_j b\in l^\infty_j L_x^\infty,\text{ and }2^j \Delta_j b\in
l^\infty_j L^1_x \, (\text{recall that } \dot W^1_1 \hookrightarrow
\dot B^{1,\infty}_1).
$$ Before applying the remaining $\partial_x$, we have a summand $P_j$
 which is such that $ 2^{\frac j 2} P_j \in l^2_j L^1_x L^2_t$ and
 frequency localized in a ball of size $2^j$, hence $\sum_j P_j \in
 \dot B^{\frac 1 2,2}_1(\mathcal{L}^2_t)$; the result follows by
 derivation. The commutator term is essentially the same, thanks to
 the following lemma.
\begin{lemme}
\label{lemme-commutation}
  Let $g(x,t)$ be such that $\|\partial_x
  g\|_{L^{p_1}_x(L^{q_\infty}_t)}<+\infty$, and $f(x,t)\in
  L^{p_\infty}_x(L^{q_2}_t)$, with $\frac 1 {p_1}+\frac 1
  {p_\infty}=1$ and $\frac 1 {q_\infty}+\frac 1 {q_2}=\frac 1 2$, then
  $h(x,t)=[\Delta_j,g] f$ is in $L^1_x(L^2_t)$.
\end{lemme}
\dem We first take $p_1=1$, $p_\infty=\infty$: set $h(x)=[\Delta_j,g]
f$, recall $\Delta_j$ is a convolution by $2^j\phi(2^j \cdot)$, and
denote $\psi(z)=z|\phi|(z)$:
\begin{equation*}\begin{aligned}
    h(x) & = \int_y 2^j\phi(2^j(x-y))(g(y)-g(x))f(y)dy \\ & = \int_{y
,\theta\in [0,1]} 2^j\phi(2^j(x-y))(x-y)g'(x+\theta(y-x))f(y) d\theta
dy \\ |h(x)| & \leq 2^{-j} \int_{y,\theta\in [0,1]}
2^j\psi(2^j(x-y))|g'(x+\theta(y-x))||f(y)| d\theta dy
\end{aligned}
\end{equation*}
and then take successively time norms and space norms,
\begin{equation*}
\begin{aligned}
 \|h(x,t)\|_{L^2_t} & \leq 2^{-j} \int_{y,\theta\in [0,1]}
   2^j\psi(2^j(x-y))\|g'(x+\theta(y-x,t))\|_{L^{q_\infty}_t}\|f(y,t)\|_{L^{q_2}_t}
   d\theta dy\\ \int_x \|h(x)\|_{L^2_t}dx & \leq 2^{-j}
   \|f\|_{L^\infty_x(L^{q_2}_t)} \int_{\smash{\theta\in [0,1],x,y}}
   2^j\psi(2^j(x-y))\|g'(x+\theta(y-x))\|_{L^{q_\infty}_t} dx dy
   d\theta\\ & \leq 2^{-j} \|f\|_{L^\infty_x(L^{q_2}_t)}
   \int_{\smash{\theta\in [0,1]},z,x} 2^j\psi(2^jz)\|g'(x+\theta
   z)\|_{L^{q_\infty}_t} dx dz d\theta\\ & \leq 2^{-j}
   \|f\|_{L^\infty_x(L^{q_2}_t)} \int_z 2^j\psi(2^jz) dz
   \|g'(x)\|_{L^1_x(L^{q_2}_t)}.
\end{aligned}
\end{equation*}
The case $p_1=\infty$, $p_\infty=1$ is identical, exchanging $f$ and
$g'$ (in fact, this would be the usual commutator estimate !). The
general case then follows by bilinear complex interpolation. \cqfd

Thus, the lemma allows us to effectively proceed with the second term
in $R_j$ as if the derivative on $\Delta_k u$ was in fact on an
$S_{k-3} b$ factor, and then it becomes a term ``like''
$$ \partial_x \sum_{k\sim j} S_k b \partial_x S_k b \Delta_k u,
$$ for which the computation done with the first term holds as
well. We are left with the third term: this is nothing but a
paraproduct
which is easily estimated: $\partial_x S_{j-3} b\in L^1_x$ and
$2^{-\frac j 2}\partial_x \Delta_j u\in L^\infty_x L^2_t$.  This
completes the proof of Proposition \ref{reminder}.\cqfd

After the paralinearization step, we perform a change of variable.  We
have, denoting by $\omega=\sqrt{\frac m 2+ (S_{j-3} b)^2}$, and $
u_j=\Delta_j u$,
$$ i\partial_t u_j+\omega(x) \partial_x ( \omega(x) \partial_x
 u_j)=R_j.
$$ Now we set $x=\phi(y)$ through $\partial_y=\omega(x)\partial_x$, in
 other words
$$ \omega(x)=\frac{ dx}{dy},\,\,\, y=\int_0^x \omega(\rho)\,
 d\rho=\phi^{-1}(x),
$$ which is a $C^1$ diffeomorphism (uniformly with respect to $j$): $
 \omega$ is bounded in the range $[\frac m 2, 2M]$. Denote by
 $v_j(y)=u_j\circ \phi(y)$ and
 $T_j(y)=R_j\circ \phi(y)$,
 \begin{equation*}
   \label{eq:v}
 i\partial_t v_j+\partial^2_y v_j= T_j(y).
 \end{equation*}
Given that our change of variable leaves $L^p$ spaces invariant, from
 Proposition \ref{reminder}, we have that
\begin{equation}\label{est.remain}
T_j \in L^1_y L^2_t, \text{ with } \|T_j\|_{L^1_y L^2_t} \lesssim
2^{\frac j 2} \mu_j,\qquad (\mu_j)_j \in l^2.
\end{equation}
By using Duhamel,
\begin{equation}
  \label{eq:duhamel-v}
 v_j=S(t)v_j(0)+\int_0^t S(t-s) T_j(y,s)ds
\end{equation}
for which we can apply Christ-Kiselev Lemma; first, let us obtain
Strichartz estimates: according to \eqref{est.remain},
\eqref{eq:duhamel-v} and Theorem~\ref{th-CK}, we obtain
\begin{equation}
  \label{eq:Striv}
\| v_j \|_{L^4_t \dot B^{\frac 1 2,2}_\infty} \lesssim
\|v_j(0)\|_{\dot H^{\frac 1 2}}+ 2^{\frac j 2} \mu_j.
\end{equation}
Now we would like to go back to $u_j$ from $v_j$.  While frequency
localizations wrt $x$ and $y$ do not commute, they ``almost'' commute.
\begin{proposition}
\label{besov-diffeo}
  Let $x=\phi(y)$ be our diffeomorphism, $|s|<1$ and $1\leq p,q\leq
  +\infty$. Then the Besov spaces $\dot B^{s,q}_p(x)$ and $\dot
  B^{s,q}_p(y)$ are identical, with equivalent norms.
\end{proposition}
\dem For any $p\in [1,+\infty]$, the $\dot W^1_p$ norms are
equivalent: the two Jacobians $|\partial_y \phi(y)|$ or $|\partial_x
\phi^{-1}(x)|$ are bounded. Therefore, with obvious notations,
$$ \|\Delta^y_j \Delta^x_k \varphi\|_p\sim 2^{-j} \|\Delta^y_j
 \Delta^x_k \varphi\|_{\dot W^1_p(y)}\lesssim 2^{-j} \|\Delta^x_k
 \varphi\|_{\dot W^1_p(x)} \lesssim 2^{k-j} \|\Delta^x_k
 \varphi\|_p\sim 2^{k-j} \|\varphi\|_p.
$$ Since $x$ and $y$ play the same part, by duality we obtain
$$ \|\Delta^y_j \Delta^x_k \varphi\|_p\lesssim 2^{-|k-j|}
 \|\varphi\|_p.
$$ This essentially allows to exchange $x$ and $y$ in Besov spaces, as
 long as we are using spaces involving strictly less than one
 derivative: say $\varphi(x) \in \dot B^{s,q}_p(x)$, then
 $\varphi(y)\in \dot B^{s,q}_p(y)$, as
\begin{equation*}
\begin{gathered}
\|\Delta^y_j \varphi\|_p  \lesssim \sum_k 2^{-|k-j|} \|\Delta^x_k
\varphi\|_p  \lesssim \sum_k 2^{-|k-j|} 2^{-sk} \varepsilon_k\\ 2^{js}
\|\Delta^y_j \varphi\|_p   \lesssim \sum_k 2^{-(1-s)|k-j|}
\varepsilon_k \lesssim  \mu_j,
\end{gathered}
\end{equation*}
where $(\mu_j )_j \in l^q$ as an $l^1-l^q$ convolution.\cqfd
\begin{rem}
Proposition \ref{besov-diffeo} is nothing but the invariance of Besov
spaces under diffeomorphism. Given that we only have a $C^1$
diffeomorphism, we are restricted to Besov spaces with $|s|<1$
regularity.
\end{rem}
Going back to \eqref{eq:Striv}, we immediately obtain by inverting the
change of variable,
$$ \|u_j \|_{L^4_t(\dot B^{\frac 1 2,2}_\infty)} \lesssim
\|u_j(0)\|_{\dot H^\frac 1 2}+ 2^{\frac j 2} \mu_j,
$$ and given that $u_j=\Delta_j u$,
$$ \|u_j \|_{L^4_t(L^\infty_x)} \lesssim \|\Delta_j u_0\|_2+ \mu_j,
$$ which, by summing over $j$, gives the desired Strichartz
estimate. All other Strichartz estimates are obtained directly in the
same way or by interpolation with the conservation of mass. This ends
the proof of Theorem \ref{Strichartz-L2}.

The strategy is exactly similar for the maximal function
estimate. Recall (see~\cite{KeRu82}) that for the (flat) Schr\"odinger
equation, we have
\begin{equation}
  \label{eq:max-flat}
\|(-\partial_y^2)^{-\frac 1 8} S(t) u_0\|_{L^4_y(L^\infty_t)}\simeq
  \left\|\int_\R e^{iy\xi-it|\xi|^2} \frac{\hat u_0}{|\xi|^{\frac 1
  4}} \,d\xi\right\|_{L^4_y(L^\infty_t)} \lesssim \|u_0\|_{L^2(\R)},
\end{equation}
from which we may obtain ( by combining \eqref{eq:max-flat} with smoothing and
Christ-Kiselev) an inhomogeneous estimate for the flat case,
$$ \|(-\partial_y^2)^{\frac 1 8}\int_0^t
S(t-s)f(s)\,ds\|_{L^4_y(L^\infty_t)} \lesssim \|f\|_{L^2_y(L^1_t)}.
$$ Therefore applying this estimate on \eqref{eq:duhamel-v} (at the
frequency-localized scale) we get
$$ \|v_j \|_{\dot B^{\frac 1 4,2}_4(L^\infty_t)} \lesssim
\|v_j(0)\|_{\dot H^\frac 1 2}+2^{\frac j 2} \mu_j,
$$ and then
$$ \|u_j \|_{L^4_x(L^\infty_t)} \lesssim 2^{\frac j 4} (\|\Delta_j
u_0\|_2+\mu_j),
$$ which we can then sum up.
\begin{rem}
  Here we are using an equivalence between Besov spaces wrt $x$ and
  Besov spaces wrt $y$ with value in $L^\infty_t$. The reader will
  easily check that the argument we used to obtain Proposition
  \ref{besov-diffeo} applies with any Besov spaces with value in
  $L^q_t$ for any $1\leq q\leq +\infty$. As an alternative, one could
  use the definition with moduli of continuity (which is the usual way
  to prove invariance by diffeomorphism) to obtain the $0<s<1$ range
  (and duality if one needs $-1<s<0$).
\end{rem}
 This completes the proof of Theorem \ref{maximal} for the special
 case $s=\frac 1 4$. We are left with shifting regularity in the
 appropriate range: but this is again nothing but a consequence of the
 equivalence from Appendix \ref{equiv-norm}. We therefore obtain the
 full range in Theorem \ref{maximal} as well as Proposition
 \ref{Strichartz-Hs}, where the restriction on $s$ follows from book
 keeping.\cqfd

\section{Application to a generalized Benjamin-Ono equation}
\label{BO}

Benjamin-Ono reads
\begin{equation}
  \label{eq:bop}
(\partial_t+H\partial_x^2) u\pm u^{p} \partial_x u=0,
\end{equation}
with real data $u_0$ at time $t=0$ (thus, it stays real). Here $H$
denotes the Hilbert transform (Fourier multiplier
$i\text{sign}(\xi)$). Given that the solution is real-valued, we can
recover it from its positive spectrum; by projecting on positive
frequencies, we get a Schr\"odinger equation. In particular,
smoothing, Strichartz, maximal function estimates are strictly the
same for both linear operators.

There are several cases of interest: mainly $p=1$, $p=2$ and $p=4$. We
will restrict ourselves to $p=4$. Other cases will be dealt with
elsewhere (\cite{BPBO}), as they require extra developments and
significantly new ideas in addition to the techniques we developed in
the present paper.

The study of the IVP for \eqref{eq:bop} with low regularity data was
initiated in \cite{KPV91,KPV94}. The best results to date were
obtained recently in \cite{MR2}, where they prove (among other results
for different $p$) \eqref{eq:bop} to be locally wellposed in
$H^{\frac 1 2^+}$. The authors were able to remove the (rather natural
with the techniques at hand) restriction on the size of the data by
adapting the renormalization procedure from \cite{TaoBO} (where global
wellposedness for the $p=1$ case is obtained in $H^1$). The same
authors proved earlier in \cite{MR1} that \eqref{eq:bop} was globally
wellposed for small data in $\dot B^{1/4,1}_4$ (and extended this
result to $\dot H^{\frac 1 4}$ in \cite{MR2}). We refer to \cite{MR2}
for a very nice presentation of the Benjamin-Ono family of equations
and of the context in which they arise.

We intend to remove the restriction on the size of the data all the
way down to $s=1/4$ (which is the scaling exponent).
\begin{theoreme}
\label{BOlocal}
  Let $u_0\in \dot H^{\frac 1 4}$, then the generalized Benjamin-Ono
  equation \eqref{eq:bop}, for $p=4$ is locally wellposed, i.e. there
  exists a time $T(u_0)$ such that a unique solution $u$ exists with
  \begin{equation*}
    \label{eq:4existence}
    u\in C_T (\dot H^{\frac 1 4})\cap \dot B^{\frac 3
    4,2}_\infty(L^2_T)\cap \dot L^4_x(L^\infty_T).
  \end{equation*}
Moreover, the flow map is locally Lipschitz.
\end{theoreme}
Combining this local wellposedness result, which is subcritical with
respect to the ``energy norm'' $\dot H^{\frac 1 2}$, with the
conservation of mass and energy,
\begin{equation*}
  \label{eq:invariant}
  \| u(t)\|_2 =\|u_0\|_2 \text{ and } E(u)=\|u\|^2_{\dot H^\frac 1
  2}\mp\frac 1 {15} \int_\R u^6=E(u_0)
\end{equation*}
and Gagliardo-Nirenberg, we also obtain global wellposedness in the
energy space when the energy controls the $\dot H^\frac 1 2$ norm,
 which occurs in the defocusing case (minus sign in \eqref{eq:bop}) or if
the $L^2$ norm is small enough (focusing: plus sign in
\eqref{eq:bop}).
\begin{theoreme}
\label{BOglobal}
  Let $u_0\in H^{\frac 1 2}$, then the defocusing generalized
  Benjamin-Ono equation (\ref{eq:bop}), for $p=4$, is globally
  wellposed, i.e. there exists a unique solution $u$ such that
  \begin{equation*}
    \label{eq:4existenceglobal}
    u\in C_T (H^\frac 1 2)\cap \dot B^{\frac 3
    4,2}_\infty(L^2_{t,\text{loc}})\cap \dot
    L^4_x(L^\infty_{t,\text{loc}}).
  \end{equation*}
\end{theoreme}
\dem We first prove Theorem \ref{BOlocal}. For local
well-posedness, the sign in \eqref{eq:bop} is irrelevant and we take
$+$ for convenience.
Let us sketch our strategy: the restriction on small data is induced
by the maximal function estimate \eqref{eq:max-flat}: even on the
linear part, $\|S(t)u_0\|_{L^4_x(L^\infty_t)}$ will be small only if
$\|u_0\|_{\dot H^\frac 1 4}$ is small as well. Here and hereafter,
$S(t)$ denote the linear operator, which we recall reduces to the
Schr\"odinger group on positive frequencies. Now, if we consider
instead the difference $S(t)u_0-u_0$, then the associated maximal
function is small provided we restrict ourselves to a small time
interval $[0,T]$:
  \begin{lemme}
\label{diffmax}
    Let $u_0\in \dot H^{\frac 1 4}$, then for any $\e>0$, there exists
    $T(u_0)$ such that
    \begin{equation}
      \label{eq:maximalT}
      \|\sup_{|t|<T} |S(t) u_0-u_0|\|_{L^4_x} < \e.
    \end{equation}
  \end{lemme}
\dem For the linear flow,
  \begin{eqnarray*}
    \|S(t)u_0-u_0\|_{L^4_x(L^\infty_T)} & \leq & \sum_{|j|<N}
    \|\Delta_j(S(t)u_0-u_0)\|_{L^4_x(L^\infty_T)}+2 \Bigl(\sum_{|j|>N}
    2^{\frac j 2} \|\Delta_j u_0\|^2_2\Bigr)^\frac 1 2\\ & \leq &
    \sum_{|j|<N} 2^{2j}\|\int_0^t S(s)\Delta_j u_0
    ds\|_{L^4_x(L^\infty_T)}+2 \Bigl(\sum_{|j|>N} 2^{\frac j 2}
    \|\Delta_j u_0\|^2_2\Bigr)^\frac 1 2\\ & \leq & T \sum_{|j|<N}
    2^{2j}\| \Delta_j u_0 \|_{L^4_x(L^\infty_T)}+2 \Bigl(\sum_{|j|>N}
    2^{\frac j 2} \|\Delta_j u_0\|^2_2\Bigr)^\frac 1 2\\ & \leq & T
    2^{2N} \| u_0 \|_{\dot H^\frac 1 4} +2 \Bigl(\sum_{|j|>N} 2^{\frac
    j 2} \|\Delta_j u_0\|^2_2\Bigr)^\frac 1 2
  \end{eqnarray*}
and by choosing first $N$ large enough and then $T$ accordingly, we
get arbitrary smallness.\cqfd

Given that local in time solutions do exist (\cite{KPV91}), we could
set up an a priori estimate and pass to the limit. However, in order
to get the flow to be Lipschitz, one has essentially to estimate
differences of solutions, and in turn this provides the required
estimates to set up a fixed point procedure.

Firstly, we proceed with an appropriate paralinearization of the
equation itself. All computations which follow are justified if we
consider smooth solutions. We have, denoting $u_j=\Delta_j u$,
$u_{\prec j}=S_{j-10} u$ and $u_{\preceq j}=S_j u$
$$ \partial_t u_j+H\Delta u_j+\Delta_j (u^4\partial_x u)=0.
$$ Rewriting $u^4\partial_x u=\partial_x(u^5)/5$ and using a
telescopic series $u=\sum_k S_k u-S_{k-1} u$, we get by standard
paraproduct-like rearrangements
\begin{align*}
 5 \Delta_j (u^4\partial_x u) = \Delta_j\partial_x(u^5) & =
   \Delta_j \Bigl((u_{\prec j})^4 \sum_{k\sim j} \partial_x
 u_k\Bigr)+\partial_x \Delta_j\Bigl(\sum_{j\lesssim k\sim k'} (u_{k'})^2
 (u_{\preceq k'})^3\Bigr)\\ & {}+\Delta_j \Bigl(\sum_{k\sim j}
 (u_{\prec j})^3 u_k \partial_x u_{\prec j}\Bigr) = \Delta_j
 \Bigl((u_{\prec j})^4 \sum_{k\sim j} \partial_x u_k\Bigr) -R_j(u).
\end{align*}
We will now consider the original equation as a system of frequency
localized equations,
\begin{equation*}
  \label{eq:microlocalbo4}
  \partial_t u_j+H\Delta u_j+ \Delta_j \Bigl(\sum_{k\sim j} (u_{\prec
  j})^4 \partial_x u_k\Bigr)=R_j(u).
\end{equation*}
 If we set $ \pi(f_1,f_2,f_3,f_4,g)= \sum_j \Delta_j \Bigl(\sum_{k\sim
  j} f_{1,\prec j}f_{2,\prec j}f_{3,\prec j}f_{4,\prec j} g_k\Bigr) $
  we can rewrite our model (abusing notations for $\pi$)
\begin{equation}
\label{eq:paraBO}
 \partial_t u+H\Delta u+\pi(u^{(4)},\partial_x u) = R(u),
\end{equation}
and we intend to solve (\ref{eq:paraBO}) by Picard iterations.

Now, let us consider $u_{L}$ the solution to the linear BO equation,
and the following linear equation:
\begin{equation*}
  \label{eq:bo4lineaire}
   \partial_t v+H\Delta v+\pi(u_{L}^{(4)},\partial_x v) = 0, \text{
 and } v_{t=0}=u_0.
\end{equation*}
At the frequency localized level, this is almost what we can handle,
except for a commutator term. Therefore we have
\begin{equation*}\begin{aligned}
    \partial_t v_j+H\Delta v_j+ (u_{L,\prec j})^4 \partial_x v_j & = 
  -\Bigl(\sum_{k\sim j}[\Delta_j, (u_{L,\prec j})^4] \partial_x
  v_k\Bigr),\\ \partial_t v_j+H\Delta v_j+ (u_{0,\prec j})^4
  \partial_x v_j & =  \Bigl((u_{0,\prec j})^4-(u_{L,\prec j})^4\Bigr)
  \partial_x v_j -\Bigl(\sum_{k\sim j}[\Delta_j, (u_{L,\prec j})^4]
  \partial_x v_k\Bigr),
\end{aligned}
\end{equation*}
for which we aim at using the estimates from Section \ref{sec1}.

The iteration map will therefore be
$$ \partial_t u_{n+1}+H\Delta u_{n+1}+\pi(u_L^{(4)},\partial_x
 u_{n+1}) = \pi(u_L^{(4)},\partial_x u_{n})-\pi(u_n^{(4)},\partial_x
 u_{n})+ R(u_n).
$$ Hence we need estimates for the linear equation
\begin{equation}
  \label{eq:bo4inhomogene}
   \partial_t v+H\Delta v+\pi(u_{L}^{(4)},\partial_x v) = f(x,t),
 \text{ and } v_{t=0}=u_0.
\end{equation}
 Restrict time to $[0,T]$ with $T$ to be chosen later, let $0^+$
 denote a small number close to $0$, and define
 $$ E_s=\cap_{0^+\leq \theta\leq 1} \BL {s+\frac{5\theta-1}4}
{\frac{4}{1-\theta}} 2 {\frac 2 \theta} \text{ as well as }
F_s=\sum_{0^+\leq \theta\leq 1,\text{finite}} \BL {s+\frac {1-3\theta} 4}
{\frac{4}{3+\theta}} 2 {\frac 2 {2-\theta}}
$$ (we left out the maximal function part, $\theta=0$ because we need
 a slightly different estimate).
\begin{proposition} \label{paraBOestimates}
Let $v$ be a solution of equation \eqref{eq:bo4inhomogene}, $u_0\in
\dot H^{s}\cap \dot H^{\frac 1 4}$ with $-3/4<s<1/2$ and $f\in
F_s$. Then there exists $T(u_0)$ such that on the time interval
$[-T,T]$, we have
\begin{equation*}
  \label{eq:EsFs}
\|v\|_{E_s} \lesssim_T \|u_0\|_{\dot H^{s}}+\| f\|_{F_s}.
\end{equation*}
Moreover,
\begin{equation*}
  \label{eq:paraBOmaxs}
\|v-u_0\|_{\BL {s-\frac 1 4} 4 2 \infty} \lesssim_T
  \|S(t)u_0-u_0\|_{\BL {s-\frac 1 4} 4 2 \infty} +\|f\|_{\BL {s-\frac
  1 2} 1 2 2}+\|u_0\|^4_4 \|v\|_{\BL {s+\frac 1 2} \infty 2 2},
\end{equation*}
and
\begin{equation*}
  \label{eq:paraBOmax}
\|v-u_0\|_{L^4_x(L^\infty_T)} \lesssim_T \|S(t)u_0-u_0\|_{L^4_x
 L^\infty_T} +\|f\|_{\BL {-\frac 1 4} 1 2 2}+\|u_0\|^4_4 \|v\|_{\BL
 {\frac 3 4} \infty 2 2}.
\end{equation*}
\end{proposition}
\dem Let us consider the equation at the frequency localized level,
$$ \partial_t v_j+H\Delta v_j+ (u_{0,\prec j})^4 \partial_x v_j =
  \Bigl((u_{0,\prec j})^4-(u_{L,\prec j})^4\Bigr) \partial_x v_j
  -\Bigl(\sum_{k\sim j}[\Delta_j, (u_{L,\prec j})^4] \partial_x
  v_k\Bigr)+f_j,
$$ and we will denote by $R_j$ the right hand side. Notice $R_j$ is
spectrally localized.  In order to connect this equation with the
model worked upon in Section \ref{sec1}, denote by
$$ b(x)=(u_{0,\prec j})^4\in L^1_x, \text{\ \ and consider \ } i
\partial_t w+ \partial_x^2 w+b(x) \partial_x w = g.
$$ By reversing the procedure we used in Section \ref{sec2}, we can
reduce the operator $\partial^2_x+b(x)\partial_x$ to $\partial_y a(y)
\partial_y$ and apply all the estimates we already know: set
$$ \frac{dy}{dx}=A(x)=\phi'(x),\text{\ \ with \
 }A(x)=\exp \bigl(\int_{-\infty}^x (u_{0,\prec j})^4(\rho)\,d\rho\bigr),
$$ then $y=\phi(x)$ is a diffeomorphism and $\sqrt a(y)=A\circ
\phi^{-1}(y)$ which insures $a\in \dot W^{1,1}$, and $a$ is
$1$-admissible. A simple calculation shows that under this change of variables,
$$ 
\partial_y a(y) \partial_y \rightarrow \partial_x^2 + b(x) \partial_x
$$
Note that everything is uniform wrt $j$. Interpolation
between all the various bounds which one can deduce from Proposition
\ref{th:smoothing-s} and Theorem \ref{maximal} yields estimates for
$w\circ\phi^{-1}$ which are identical to the flat case (or, to get a
better sense of perspective, to linear estimates for the linear
Benjamin-Ono equation, see e.g. \cite{MR1}):
$$ \|w\circ\phi^{-1}\|_{E_s} \lesssim \|w_0\circ\phi^{-1}\|_{\dot
H^{s}}+\| g\circ\phi^{-1}\|_{F_s},
$$ with $-3/4<s<1/2$. Using Proposition \ref{besov-diffeo}, we can
revert back to the $x$ variable and obtain the {\em exact} same
estimates for $ w$:
$$ \|w\|_{E_s} \lesssim \|w_0\|_{\dot H^{s}}+\| g\|_{F_s}.
$$ Recalling that $w=v_j=\Delta_j v$ and $g=R_j$ is frequency
localized as well, hence for any $0^+\leq \theta\leq 1$,
$$ 2^{j(s+\frac{5\theta-1}4)}
\|v_j\|_{L^{\frac{4}{1-\theta}}_x(L^{\frac 2 \theta}_t)} \lesssim
2^{js}\|u_{0,j}\|_2+\sum_{k,\text{finite}} 2^{j( {s+\frac
{1-3\theta_k} 4} )} \|R_j\|_{L^{\frac{4}{3+\theta_k}}_x(L^{\frac 2
{2-\theta_k}}_t)}.
$$ All is left is to estimate $R_j$ in order to contract the $v_j$
term:
$$ R_j= \Bigl((u_{0,\prec j})^4-(u_{L,\prec j})^4\Bigr) \partial_x v_j
-\Bigl(\sum_{k\sim j}[\Delta_j, (u_{L,\prec j})^4] \partial_x
v_k\Bigr)+f_j.
$$ From the smoothing estimate for the flat Schr\"odinger equation and
Lemma \ref{diffmax}, there exist $T(u_0)$ such that
\begin{equation}
  \label{eq:petit}
\Bigl( \sum_j (2^{\frac {-j}{4}}\|\partial_x u_{L,\preceq
    j}\|_{L^\infty_x(L^2_T)})^2\Bigr)^\frac 1 2 +\|u_{L,\prec
    j}-u_{0,\prec j}\|_{L^4_x(L^\infty_T)}) < {\eta(u_0)},
\end{equation}
where $\eta(u_0)$ can be made as small as needed by choice of a
smaller $T(u_0)$.  This allows to write, picking a $\theta_k$ close to
$1$ and abusing notations,
\begin{align*}
  2^{js}\| v_j\|_{L^\infty_x(L^2_T)})+2^{(s^-)j}\|
v_j&\|_{L^{\infty^-}_x(L^{2^+}_T)} \leq \frac 1 2 2^{(s-1)j}\|
\partial_x v_j\|_{L^\infty_x(L^2_T)}\\ &{}+\frac 1 {2K}
\sum_{|l-j|<K}2^{(s^-)l}\| w_l\|_{L^{\infty^-}_x(L^{2^+}_T)}+
2^{(s-1)j}\|f_j\|_{L^1_x(L^2_T)},
\end{align*}
where we used Lemma \ref{lemme-commutation} to estimate the commutator
with $2^{-(1^-)j} \partial_x u_{L,\prec j} \in
L^{4^+}_x(L^{\infty^-}_T)$ small enough by \eqref{eq:petit} and
interpolation with $u_L\in L^4_x(L^\infty_T)$.  We have therefore
obtained, after summing over $j$,
$$ \| v\|_{E_s}\leq C(u_0)( \|u_0\|_{\dot H^{s}_2}+\|f\|_{F_s}).
$$ We only have a local in time estimate for the linearized equation,
but it depends only on the data and nothing else, through lemma
\ref{diffmax}. At our desired level of regularity, namely $s=3/4$,
$$ \|v\|_{\BL {\frac 3 4} \infty 2 2}\leq C(u_0)(\|f\|_{\BL {-\frac 1
 4} 1 2 2}+\|u_{0}\|_{\dot H^{\frac 1 2}}).
$$ We also need the maximal function, or more accurately, $v-u_0$: but
 this is now very easy, simply reverting back to writing ($S(t)$ being 
here the group associated to the linear BO)
$$ v=u_L+\int_0^t S(t-s) (f-\pi(u_L,\partial_x v)) ds,
$$ and we therefore get (using the third case in Theorem \ref{th-CK}
for the special case $s=0$)
$$ \|v-u_0\|_{\BL {s-\frac 1 4} 4 2 \infty}\lesssim \|u_L-u_0\|_{\BL
{s-\frac 1 4} 4 2 \infty}+ \|f\|_{\BL {s-\frac 1 2} 1 2
2}+\|u_0\|^4_{4} \|v\|_{ \BL {s-\frac 1 4} \infty 2 2},
$$ and
$$ \|v-u_0\|_{L^4_x(L^\infty_T)}\lesssim \|u_L-u_0\|_{L^4_x(L^
\infty_T)}+ \|f\|_{\BL {-\frac 1 4} 1 2 2}+\|u_0\|^4_{4} \|v\|_{ \BL
{\frac 3 4} \infty 2 2}.
$$ This achieves the proof of Proposition \ref{paraBOestimates}.\cqfd

Everything is now ready for a contraction in a complete metric space,
which will be the intersection of two balls,
$$ B_M(u_0,T)=\{u \text{ s.t. } \|u-u_0\|_{\dot
B^{0,1}_4(L^\infty_T)}<\varepsilon(u_0)\},
$$ and
$$ B_S(u_0,T)=\{u \text{ s.t. } \|u\|_{\dot B^{\frac 3
4,1}_\infty(L^2_T)}<\varepsilon(u_0)\}.
$$ We first check that the mapping $K$ is from $B_M\cap B_S$ to
itself, where $K(v)=u$ with
$$ \partial_t u+H\Delta u+\pi(u_L^{(4)},\partial_x u) =
 \pi(u_L^{(4)},\partial_x v)-\pi(v^{(4)},\partial_x v)+ R(v).
$$ For this we use Proposition \ref{paraBOestimates} with $s=3/4$ and
standard (para)product estimates. The $B_S$ part is trivial (one
doesn't even need to take advantage of the difference on the
right). The $B_M$ part follows from the ability to factor an $u_L-u$
while rewriting the difference of the $\pi$ on the right.

The next step is then to contract, i.e. estimate $K(v_1)-K(v_2)$ in
terms of $v_1-v_2$. But this is again trivial given we have a
multilinear operator, it will be exactly as the $v\ra u$
mapping. This ends the proof of Theorem \ref{BOlocal}.\cqfd

We now briefly sketch the proof of Theorem \ref{BOglobal}. We now have
a minus sign in $\eqref{eq:bop}$ but this doesn't change the local in
time contraction. Given a
datum in the (inhomogeneous) space $H^s$, with $s>1/4$, a standard
modification of the fixed point provides that the solution $u$
is $C_t(H^s)$. In order to iterate whenever $s=1/2$, we need to check
that the local time $T(u_0)$ can be repeatedly chosen in a uniform
way. All is required is an appropriate modification of Lemma
\ref{diffmax}: recall we can write
$$ \sum_j \|\Delta_j (S(t)u_0-u_0)\|_{L^4_x(L^\infty_T)} \leq T 2^{2N}
\| u_0 \|_{\dot H^{\frac 1 4}} +2 \Bigl(\sum_{|j|>N} 2^{\frac j 4}
\|\Delta_j u_0\|_2\Bigr),
$$ from which we get, taking advantage of $u_0\in L^2\cap \dot
H^{\frac 1 2}$,
$$ \sum_j \|\Delta_j (S(t)u_0-u_0)\|_{L^4_x(L^\infty_T)} \leq T 2^{2N}
\Bigl(\|u_0\|_2 \| u_0 \|_{\dot H^{\frac 1 2}}\Bigr)^{\frac 1 2}+
2^{-\frac N 4} (\|u_0\|_2+ \| u_0 \|_{\dot H^{\frac 1 2}}).
$$ Obviously, picking $T=2^{-\frac 9 4 N}$ gives the bound $ 2^{-\frac
N 4} (\|u_0\|_2+ \| u_0 \|_{\dot H^{\frac 1 2}})$, which by an
appropriate choice of $N$ can be made as small as we need with respect
to $(\|u_0\|_2+ \| u_0 \|_{\dot H^{\frac 1 2}})$. However, both the
$L^2$ and $\dot H^\frac 1 2$ norms are controlled, thus the local time
$T(u_0)$ is uniform and we can iterate the local existence result to a
global result.\cqfd

\appendix
\section{Localization with respect to $\partial_x$ versus localization
  with respect to $(-\partial_x(a(x)\partial_x))^\frac 1 2$}
\label{equiv-norm}

\subsection{The heat flow associated with
  $-\partial_x(a(x)\partial_x$}
\label{sec:heat}
We would like to define an analog of the Littlewood-Paley operator
$\Delta_j$, but using $A=-\partial_x(a(x)\partial_x$ rather than
$-\partial^2_x$. In the first 2 sections, this turns out to be useful
because such a localization wrt $A$ will commute with the
Schr\"odinger flow. Through spectral calculus, we can easily define
$\phi(A)$ for a smooth $\phi$, but we need various properties on $L^p$
spaces for all $1\leq p\leq +\infty$, which requires a bit more of
real analysis. Fortunately, all the results we need are more or less
direct consequences of (part of) earlier work related to the Kato
conjecture, and we simply give a short recollection of the main facts
we need, skipping details and referring to \cite{AT,AMT}. We call
$S_A(t)$ the heat flow, namely $S_A(t) f$ solves
\begin{equation}
  \label{eq:heatflow}
  \partial_t g+A g=0,\text{ with } g(0)=f,
\end{equation}
and define $\Delta^A_j f=4^{-j} AS_A(4^{-j}) f$. Again, in $L^2$ all
of this makes sense through spectral considerations, and were $a$ to
be just $1$, we would just get a localization operator based on the
Mexican hat $\xi^2\exp -\xi^2$. In \cite{AT}, such a semi-group
$S_A(t)$ is proved to be analytic, and moreover the square-root of $A$
can be factorized as $R\partial_x$, where $R$ is a Calderon-Zygmund
operator, under rather mild hypothesis: $a\in L^\infty$, complex
valued, with $\Re a>1$. On the other hand, in \cite{AMT}, the authors
prove Gaussian bounds for the kernel of the semi-group as well as its
derivatives, and this provides everything which is needed here. Such
bounds are obtained through the following strategy:
\begin{itemize}
\item Derive bounds for the operator $(1+A)^{-1}$: given that it maps
  $H^{-1}$ to $H^1$, it follows that it maps $L^1$ to $L^\infty$ by
  Sobolev embeddings.
\item Obtain bounds for $(\lambda+A)^{-1}$, $\Re \lambda>0$, by
  rescaling, given the hypothesis on $a$ are invariant.
\item Obtain bounds for $A (1+A)^{-1}$ by algebraic manipulations,
  proving it maps $L^1$ to $L^\infty$.
\item Obtain again an $L^1-L^\infty$ bound for $\partial_x (1+A)^{-1}$
by ``interpolation'' between the two previous bounds. This specific
bound we did prove directly in Section \ref{sec1}, namely
\eqref{eq:resestbis}.
\item Use a nifty trick (see Davies (\cite{Da1})): remark that
  provided $\omega$ is sufficiently small (wrt the lower bound of $\Re
  a$), all previous estimates hold as well for
  $$A_\omega=\exp(\omega \cdot) A \exp(-\omega \cdot).$$ Then any of
 the new kernels $K_\theta(x,y)$ are just $K(x,y)\exp(-\omega |x-y|)$,
 which gives exponential decay pointwise from the $L^1-L^\infty$
 bound.
\item Use the representation of $S_A(t)$ in term of $R_\lambda
  (A)=(\lambda+A)^{-1}$ to obtain that $S_A(t)$ maps $L^1$ to
  $L^\infty$ and that its kernel verifies Gaussian bounds, as well as
  its derivatives.
\end{itemize}
We can summarize with the following proposition.
\begin{proposition}[\cite{AMT}]
  Let $K_A(x,y,t)$ be the kernel of the heat flow $S_A(t)$.  There
  exists $c$ depending only on the lower bound of $\Re a$ and its
  $L^\infty$ norm, such that
  \begin{equation}
    \label{eq:K1}
    |K_A(x,y,t)|\lesssim \frac 1 {\sqrt t} e^{c\frac{-|x-y|^2} t},
  \end{equation}

 \begin{equation}
    \label{eq:K2}
|\partial_y K_A(x,y,t)|+ |\partial_x K_A(x,y,t)|\lesssim \frac 1 {t}
 e^{c\frac{-|x-y|^2} t},
  \end{equation}
and
 \begin{equation}
    \label{eq:K3}
    | A K_A(x,y,t)|\lesssim \frac 1 {t^{\frac 3 2}}
      e^{c\frac{-|x-y|^2} t}.
  \end{equation}
\end{proposition}
Once we have all the Gaussian bounds, it becomes very easy to prove
that $S_A(t)$ is continuous on $L^p$ (from \eqref{eq:K1}), as well as
$\Delta^A_j$ (from \eqref{eq:K3}). We are, in effect, reduced to the
usual heat equation, with appropriate Bernstein type inequalities.
\subsection{Equivalence of Besov norms}
\label{sec:besov-eq}
We first define Besov spaces using the $A$ localization rather the
usual one:
\begin{definition}
\label{dA}
Let $f$ be in $\mathcal{S}'(\mathbb{ R}^{n})$, $s<1$. We say $f$
belongs to $\dot B^{s,q}_{p,A}$ if and only if
\begin{itemize}
\item The partial sum $ \sum^{m}_{-m} \Delta^A_{j}(f)$ converges to
$f$ as a tempered distribution (modulo constants if $s\geq 1/p, q>
1$).
\item The sequence $\varepsilon_{j} = 2^{js}\| \Delta^A_{j}
(f)\|_{L^{p}}$ belongs to $l^{q}$.
\end{itemize}  
\end{definition}
Alternatively, one could replace the discrete sum with a continuous
one, which is somewhat more appropriate when using the heat flow. Both
can be proved to be equivalent, exactly as in the usual situation.

Now, our aim is to prove these spaces to be equivalent to the ones
defined by Definition \ref{d1}. In order to achieve this, we would
like to estimate $\Pi_{jk}=\Delta^A_j \Delta_k$ and its adjoint. The
adjoint can be dealt with by duality, so we focus on $\Pi_{jk}$: there
are obviously 2 cases,
\begin{itemize}
\item when $j>k$, we write
$$ \Pi_{jk}=4^{-j} S_A(4^{-j}) \partial_x a(x) \partial_x \Delta_k,
$$ which immediately yields, for any $1\leq p\leq +\infty$,
$$
  \begin{aligned}
  \|\Pi_{jk} f\|_p & = 2^{-j} \| S_A(4^{-j}) 2^{-j} \partial_x a(x)
  \partial_x \Delta_k f\|_p \\ & \lesssim 2^{-j}\| a(x) \partial_x
  \Delta_k f\|_p\lesssim 2^{-j}\| \partial_x \Delta_k f\|_p \lesssim
  2^{k-j} \|\Delta_k f\|_p,
  \end{aligned}
$$ where we used the bound \eqref{eq:K2} on $S_A(1) \partial_x$.
\item In the same spirit, when $k>j$,
$$ \Pi_{jk}=4^{-j} S_A\left(\frac{4^{-j}}{2}\right)
S_A\left(\frac{4^{-j}} 2\right) \partial_x (\partial_x)^{-1} \Delta_k,
$$ and then
$$ \|\Pi_{jk} f\|_p \lesssim 2^{j} \| S_A(\frac{4^{-j}} 2) 2^{-j}
  \partial_x (\partial_x)^{-1} \Delta_k f\|_p \lesssim 2^{j}\|
  (\partial_x)^{-1} \Delta_k f\|_p\lesssim 2^{j-k}\|\Delta_k f\|_p,
$$ where we used (again) the bound \eqref{eq:K2} on $S_A(1)
\partial_x$.
\end{itemize}
Therefore,
\begin{proposition}
  Let $|s|<1$, $1\leq p,q\leq +\infty$, then $\dot B^{s,q}_p$ and
  $\dot B^{s,q}_{p,A}$ are identical, with equivalence of norms.
\end{proposition}
\begin{rem}
In previous sections, we actually used Besov spaces taking values in the
  separable Hilbert space $L^2_t$: as a matter of fact, one can reduce
  to the scalar case by projecting over an Hilbert basis, hence the
  Hilbert-valued result holds as well.
\end{rem}
\section{Christ-Kiselev lemma for reversed norms}
\label{Christ-Kiselev}
\begin{figure}[ht]
$$\ecriture{\includegraphics[width=6cm]{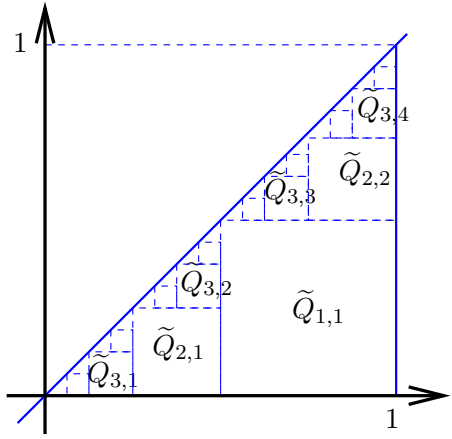}}{
\aat{42}{1}{$1$}\aat{2}{42}{$1$}\aat{32}{13}{$\widetilde{Q}_{1,1}$}\aat{17}{9}{$\widetilde{Q}_{2,1}$}\aat{37}{28}{$\widetilde{Q}_{2,2}$}\aat{10}{6}{$\widetilde{Q}_{3,1}$}\aat{20}{16}{$\widetilde{Q}_{3,2}$}\aat{29}{26}{$\widetilde{Q}_{3,3}$}\aat{39}{35}{$\widetilde{Q}_{3,4}$}
}$$\caption{Decomposition of a triangle as a union of squares}
\label{fig1}
\end{figure}

As observed in~\cite{P2nls} and further exploited in \cite{MR1},
Christ and Kiselev Lemma works also with reversed norms. In this
appendix, we prove the versions of this result we need in the previous
sections. The proof is very much inspired from~\cite{Christ-Kiselev}:
\begin{theoreme}
\label{th-CK}
\begin{itemize}
\item Let $1\leq \max(p,q)<r\leq +\infty$, $B$ a Banach space, and $T$
a bounded operator from $L^p(\mathbb{R}_y; L^q( \mathbb{R}_s))$ to
$L^r( \mathbb{R}_t; B)$ with norm $C$. Let $K(y,s,t)$ be its kernel,
and $K\in L^1_{\text{loc}}( \mathbb{R}^3_{y,s,t})$ taking values in
the class of bounded operators on $B$.  Define $T_R$ to be the
operator with kernel $1_{s<t} K(y,s,t)$. Then $T_R$ is bounded from
$L^p(\mathbb{R}_y; L^q( \mathbb{R}_s))$ to $L^r( \mathbb{R}_t;B)$ with
norm smaller than $C/(1- 2^{1/r-1/\max(p,q)})$.
\item If $\max(p,q) < \min(\alpha, \beta)$ and $T$ is a bounded
operator from the space $L^p(\mathbb{R}_y; L^q( \mathbb{R}_s))$ to
$L^\alpha( \mathbb{R}_x; L^\beta(\R_t))$ with norm $C$. Let $K(y,s,x,t)$
be its kernel, and $K\in L^1_{\text{loc}}( \mathbb{R}^3_{y,s,x,t})$.
Define $T_R$ to be the operator with kernel $1_{s<t} K(y,s,x,t)$. Then
$T_R$ is bounded from $L^p(\mathbb{R}_y; L^q( \mathbb{R}_s))$ to $L^\alpha(
\mathbb{R}_x;L^\beta(\R_t))$ with norm smaller than ${C}/({1- 2^{1/\min(\alpha,
\beta)-1/\max(p,q)}})$.
\item If $T$ is a bounded operator from $\dot{B}^{0,2}_1( \mathcal{L}^2_t)$ to
$L^4( \mathbb{R}_x; L^\infty(\R_t))$ with norm $C$. Let $K(y,s,x,t)$ be
its kernel, and $K\in L^1_{\text{loc}}( \mathbb{R}^3_{y,s,x,t})$.
Define $T_R$ to be the operator with kernel $1_{s<t} K(y,s,x,t)$. Then
$T_R$ is bounded from $\dot{B}^{0,2}_1( L^2_t)$ to $L^4( \mathbb{R}_x;
L^\infty(\R_t))$ with norm smaller than ${C}/({1- 2^{-1/4}})$.
\end{itemize}
\end{theoreme}
\dem We study the first case in Theorem~\ref{th-CK}. For
any (smooth) function $f \in L^{p}( \mathbb{R}_y;L^q(
\mathbb{R}_s))$ such that $\|f\|_{L^{p}(
\mathbb{R}_y;L^q(\mathbb{R}_s))}=1$, the function $F(t)= \|1_{s<t}f(s,y)\|^p_{L^{p}( \mathbb{R}_y;L^q(\mathbb{R}_s))}$ is an increasing function from $\mathbb{R}$ to $[0,1]$, and without
loss of generality we can take it to be injective (hence,
invertible). We have
\begin{lemme}\label{lem.decomp}
 For any $f \in L^{p}( \mathbb{R}_y;L^q(\mathbb{R}_s))$, such that
$\|f\|_{L^{p}( \mathbb{R}_y;L^q(\mathbb{R}_s))}=1$,
\begin{equation}
\|1_{F^{-1}(]a,b[)}f\|_{L^{p}( \mathbb{R}_y;L^q(\mathbb{R}_s))}\leq C
|b-a|^{\frac 1 {\max(p,q)}}
\end{equation}
\end{lemme}
Indeed denote by $]t_a, t_b[= F^{-1} (]a,b[)$ and
$$G(t,x)= \left(\int_{s<t} |f(s,x)|^q ds\right)^{\frac 1 q} $$
\begin{enumerate}
\item If $p\geq q$, using that for $a,b\geq 0$ we have
$(a+b)^{p/q}\geq a^{p/q} + b^{p/q}$ we obtain
\begin{equation*}
\begin{aligned}
\|1_{F^{-1}(]a,b[)}f\|^p_{L^{p}(\mathbb{R}_y;L^q(\mathbb{R}_s))}&=
\int_x \left( \int_{t_a\leq s \leq t_b} |f(s,x)|^q ds \right)^{\frac p
q}dx\\ & = \int_x \left( \int_{ s \leq t_b} |f(s,x)|^q ds - \int_{ s
\leq t_a} |f(s,x)|^q ds \right)^{\frac p q}dx\\ & \leq \int_x \left(
\int_{ s \leq t_b} |f(s,x)|^q ds\right)^{\frac p q}- \left(\int_{ s
\leq t_a} |f(s,x)|^q ds \right)^{\frac p q}dx\\ &\leq F(t_b)- F(t_a)=
b-a
\end{aligned}
\end{equation*}
\item If $p\leq q$, using that for $x,y\geq 0$, $(x^q-y^q) \leq \frac
q p (x^p -y^p)(\max(x,y)^{q-p}$, we obtain
\begin{multline*}
\|1_{F^{-1}(]a,b[)}f\|^p_{L^{p}(\mathbb{R}_y;L^q(\mathbb{R}_s))}= \int_x \left( \int_{t_a\leq s \leq t_b} |f(s,x)|^q ds
\right)^{\frac p q}dx\\ 
\begin{aligned}
  & = \int_x \left( G(t_b,x)^q- G(t_a,x)^q \right)^{\frac p
q}dx\\ &  \leq \frac p q \int_x \left( G(t_b,x)^p- G(t_a,x)^p
\right)^{\frac p q}\left( G(t_b,x)^{q-p} \right)^{\frac p q}dx\\ 
 & \leq \frac p q \left(\int_x G(t_b,x)^p- G(t_a,x)^p dx \right)^{\frac p
q}\left(\int_x G(t_b,x)^{(q-p)\frac q {(q-p)}}dx \right)^{\frac {q-p}
q}\\
& \leq \frac p q \left(b-a \right)^{\frac p q}\|f\|^{\frac {q-p}
q}_{L^{p}( \mathbb{R}_y;L^q(\mathbb{R}_s))} \leq \frac p q \left(b-a
\right)^{\frac p q}.
\end{aligned}
\end{multline*}
\end{enumerate}
Consider now the dyadic decomposition of the real axis given by
$$\mathbb{R}= ]- \infty, t_{n,1}] \cup ]t_{n,1}, t_{n,2}]\cup \dots
\cup ]t_{n, 2^n-1}, + \infty[= \cup_{j=1}^{2^n} I_j$$ such that
$$\|f\|^{r'}_{L^{r'}( ]t_{n,j}, t_{n, j+1}];B)} = 2^{-n}$$ with the
convention $t_{n,0}= - \infty$ and $t_{n,2^{n+1}}= + \infty$. Remark
that $F(t_{n,j})=j2^{-n}$ is the usual dyadic decomposition of the
interval $[0,1[$.  We have
$$1 _{s<t}= \sum_{n=1}^{+\infty} \sum_{j=1}^{2^{n-1}}1_{(s,t)\in
Q_{n,j}}$$ where $(s,t) \in Q_{n,j} \Leftrightarrow (F(s), F(t) )\in
\widetilde{Q}_{n,j}$ and $\widetilde{Q}_{n,j}$ is as in
Figure~\ref{fig1}.

 Remark that $1_{(s,t)\in Q_{n,j}}= 1_{t\in I_{n,j}} 1_{s\in I'_{n,j}}$ for
suitable dyadic intervals $I_{n,j}$ and $I'_{n,j}$.

We are now ready to prove the main estimate:
\begin{equation*}
\|T_R f\|_{L^{r}(\mathbb{R}_t;B)}=
\|\sum_n \sum_{j=1}^{2^{n-1}} T_{n,j}f\|_{L^{r}(\mathbb{R}_t;
B)}
\end{equation*}
 where the kernel of the operator $T_{n,j} $ is equal
to $K(y,s,t)\times 1_{(s,t)\in Q_{n,j}}$. Consequently $T_{n,j}$ is
(uniformly) bounded from $L^{r'}( \mathbb{R}_t;B)$ to
$L^{p'}(\mathbb{R}_y; L^{q'}( \mathbb{R}_s))$ with norm smaller
than~$C$.

Since ${p'}\geq {q'}$ and for fixed $n$, the functions $T_{n,j}f$ have
disjoint support (in the variable $t$) we have
\begin{equation*}\begin{aligned}
\|T_R f\|_{L^{r}(\mathbb{R}_t;B)}&\leq C \sum_n
\Bigl(\sum_{j=1}^{2^{n-1}}\|1_{s\in ]t_{n,j},
t_{n,j+1}[}f\|^r_{L^{p}(\mathbb{R}_y;L^q(\mathbb{R}_s))}\Bigr)^{\frac
1 r}\\ &\leq C \sum_n \Bigl(\sum_{j=1}^{2^{n-1}}2^{-\frac{nr}
{\max(p,q)}}\Bigr)^{\frac1 r} =(1- 2^{\frac 1 r - \frac 1
{\max(p,q)}})^{-1}.
\end{aligned}
\end{equation*}
We now study the second case in Theorem~\ref{Christ-Kiselev}. The
proof relies on
\begin{lemme}\label{lem.base}
Assume that $(f_k)_{k\in \mathbb{N}}$ have disjoint supports in
$t$. Then
$$\| \sum_k f_k \|_{L^\alpha(\R_x;L^\beta(\R_t))} \leq \Bigl( \sum_k \|f_k
\|^{\min ( \alpha, \beta)}_{L^\alpha(\R_x;L^\beta(\R_t))}\Bigr) ^{\frac 1
{\min ( \alpha, \beta)}}
$$
\end{lemme}
We distinguish two cases:
\begin{itemize}
\item{$\beta \geq \alpha$}
\begin{equation*}
\| \sum_k f_k \|_{L^\alpha(\R_x;L^\beta(\R_t))} = \Bigl(\int_x \Bigl( \sum_k
\int_t |f_k |(t,x)^\beta dt \Bigr)^{\alpha/\beta}dx \Bigr)^{1/\alpha}
\end{equation*} 
but since $\alpha \leq \beta$, we have $(\sum_ka_k)^{\alpha /\beta}
\leq \sum_k a_k ^{\alpha /\beta}$ and we obtain
\begin{equation*}
\| \sum_k f_k \|_{L^\alpha(\R_x;L^\beta(\R_t))} \leq \Bigl(\int_x \Bigl(
\sum_k \int_t |f_k |(t,x)^\beta dt \Bigr)^{\alpha/\beta}dx
\Bigr)^{1/\alpha}
\end{equation*} 
\item{$\beta \leq \alpha$}
\begin{multline*}
\| \sum_k f_k \|_{L^\alpha(\R_x;L^\beta(\R_t))} = \Bigl\|\Bigl(\sum_k \int_t
|f_k|^\beta\Bigr)\|_{L_x^{\alpha/ \beta}}^{1/\beta} \leq \Bigl(\sum_k
\|\Bigl( \int_t |f_k|^\beta\Bigr)\|_{L_x^{\alpha/ \beta}}\Bigr)^{\frac
1 \beta}\\ \leq \Bigl( \sum_k \|f_k \|^{ \beta}_{L^\alpha(\R_x;
L^\beta(\R_t))}\Bigr) ^{\frac 1 {\beta}}
\end{multline*}
\end{itemize}
To prove the second case in Theorem~\ref{Christ-Kiselev}, we use the
same dyadic decomposition of $\mathbb{R}$ as before and use
Lemma~\ref{lem.base} to estimate $\|T_R f\|_{L^\alpha_x,
L^\beta_t}$. This gives
$$\|T_R f\|_{L^\alpha(\R_x; L^\beta(\R_t))}\leq \sum_n
\Bigl(\sum_{j=1}^{2^{n-1}}\|T_{n,j}f\|^{\min(\alpha,
\beta)}_{L^{\alpha}(\R_x;L^\beta(\R_t))}\Bigr)^{\frac 1
{\min(\alpha, \beta)}}
$$ and we conclude as in the previous case.

Finally, to prove the last case in Theorem~\ref{Christ-Kiselev}, we
need to combine Lemma~\ref{lem.base} with $\alpha =4, \beta= + \infty$
to deal with the $L^4(\R_x; L^\infty(\R_t))$ norm with a choice of a suitable
dyadic decomposition and prove the analog of Lemma~\ref{lem.decomp}
for the Besov space $\dot B^{0,2}_1(\R_x)$.  The dyadic decomposition is based on
$$F(t)= \sum_j \|\Bigl( \int_{s<t} |\Delta_j f(s) |^2 ds
\Bigr)^{1/2}\|_{L^1_x}^2 = \sum_j \gamma_j(t) ^2.
$$
\begin{lemme} For any function $f$ such that $\|f\|_{B^{0,2}_1(\mathcal{L}^2_t)} = \Bigl(\sum_j\| \Delta_j
f\|^2_{L^1_x; L^2_t}\Bigr) ^{1/2} =1$ we have $\|1_{F^{-1} (]a,b[)}f\|_{B^{0,2}_1(\mathcal{L}^2_t)}\leq C
(b-a)^{\frac 1 2}$.
\end{lemme}
\dem Denote by $J_j(t,x)=(\int_{s<t} |\Delta_j
f(s)|^{2}\,ds)^{\frac{1}{2}}$. Then (using $2\geq 1$)
\begin{multline}
  \|\Delta_j \chi_{F^{-1}(I)}(s)f(s)\|^{2}_{L^{2}_t}  =
  \int_{t_a}^{t_b} |\Delta_j f(s)|^{2}\,ds  = J_j(t_b,x)^{2}-
  J_j(t_a,x)^{2},\\  \leq (J_j(t_b,x)-
  J_j(t_a,x))(J_j(t_b,x)+J_j(t_a,x)).
\end{multline}

Then we add the $L^{1}_x$ norm, to get (using Cauchy-Schwarz at the
second line)
\begin{multline}
  \int_x \|\Delta_j \chi_{F^{-1}(I)}(s)f(s)\|_{L^{2}_t}\,dx \lesssim
   \int_x (J_j(t_b,x)- J_j(t_a,x))^{\frac{1}{2}}
   (J_j(t_b,x)+J_j(t_a,x))^{\frac{1}{2}} \,dx\\ \lesssim \Bigl(\int_x
   J_j(t_b,x)- J_j(t_a,x)\,dx\Bigr)^{\frac{1}{2}}\Bigl( \int_x
   (J_j(t_b,x)+J_j(t_a,x))\,dx\Bigr)^\frac 1 2
\end{multline}
and consequently
\begin{multline}
 \sum_j \Bigl(\int_x \|\Delta_j
  \chi_{F^{-1}(I)}(s)f(s)\|_{L^{2}_t}\,dx \Bigr)^2 \lesssim \sum_j
  (\gamma_j(t_b)-\gamma_j(t_a))(\gamma_j (t_b)+\gamma_j(t_a))\\
  \lesssim \sum_j (\gamma^2_j(t_b)-\gamma^2_j(t_a))=F(t_b)-F(t_a)=|I|.
\end{multline}
The rest of the proof of Theorem~\ref{Christ-Kiselev} is as in the
previous cases.

\section{A singular metric}
\label{sec.A2}
\begin{figure}[h]
$$\ecriture{\includegraphics[width=5cm]{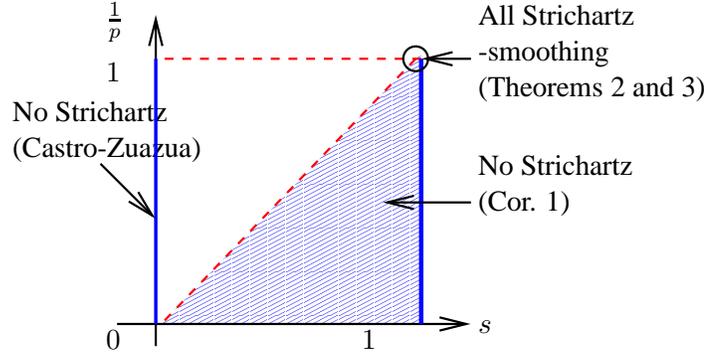}}{
\aat{35}{0}{$1$}\aat{2}{0}{$0$}\aat{2}{35}{$1$}\aat{50}{2}{$s$}\aat{2}{42}{$\frac
1 p$}\aat{50}{33}{\vbox{\hbox{All Strichartz} \hbox{-smoothing}
\hbox{(Theorems~\ref{Strichartz-L2}
and~\ref{maximal})}}}\aat{-10}{25}{\vbox{\hbox{No Strichartz}
\hbox{(Castro-Zuazua)}}}\aat{50}{18}{\vbox{\hbox{No
Strichartz}\hbox{(Cor. \ref{cor.nostrich})}}}}$$
\caption{Range of regularity $W^{s,p}$}
\label{fig.2}
\end{figure}

In this section we construct a metric on $\mathbb{R}$, which is in
$W^{s,1}$ for any $0\leq s <1$ (but not in BV), bounded from below and
above and for which no smoothing estimate and no (non trivial)
Strichartz estimates hold. In fact this construction is a
simplification of an argument of Castro and Zuazua~\cite{CaZu} (whose
proof relies in turn upon some related works in semi-classical
analysis and unique continuation theories), who, in the context of
wave equations, provide counter examples with $C^{0, \alpha}, 0\leq
\alpha <1$ metrics (continuous H\"older of exponent $\alpha$
metrics). As noticed by Castro and Zuazua, these counter examples
extend to our setting. Figure~\ref{fig.2} shows the range where full
Strichartz/smoothing are true or no Strichartz/smoothing holds. A most
interesting range of regularity is $a\in W^{s, \frac 1 s}$ and in
particular $H^{1/2} = W^{1/2, 2}$ because these regularities are scale
invariant. A natural question would be to ask whether some
Strichartz/smoothing estimates might hold (possibly with derivatives
loss) at these levels of regularity. Remark that neither our counter
examples nor Castro-Zuazua's lie in this range (except for $s=0$).
\begin{proposition}\label{prop.quasi}
There exist a metric $\beta(x)\in W^{s,1}$ for any $0\leq s <1$, bounded from below and above $0<m\leq \beta(x) \leq M$ (so that $\beta \in W^{s,p}, s< 1/p$), a
sequence of functions $\phi_k\in C^\infty_0(]2^{-k- 1/2}, 2^{-k+
1/2}[)$ and a sequence $( x_k = 2^{-k}, \lambda_k= 2^k k)$ such that
 \begin{equation}\label{eq.quasi}
(\partial_x \beta(x) \partial_x + \lambda_k) \phi_k = \mathcal{O} (
\lambda_k ^{- \infty})_{H^1}
\end{equation}
 $$\|\phi_k\|_{L^2} =1,$$
\end{proposition}
\begin{cor}\label{cor.nostrich}
For the density constructed above, we have for any $r<(q-2)/2q$
(recall that by the usual Sobolev embedding, $H^{(q-2)/2q} \rightarrow
L^q$),
\begin{equation}\label{eq.nostrich}
\lim_{k\rightarrow + \infty} \frac {\| e^{ it (\partial_x \beta(x)
\partial_x)} \phi_k \|_{L^1(- \varepsilon, \varepsilon); L^q(
\mathbb{R})}} { \|\phi_k \|_{H^r}} = + \infty
\end{equation}
\end{cor}
We first show that Proposition~\ref{prop.quasi}
implies~\eqref{eq.nostrich}.  According to~\eqref{eq.quasi},
$\|\phi_k\|_{H^1} \leq C \lambda _k$ and, by
interpolation,
\begin{equation}\label{eq.sobo}
\|\phi_k\|_{H^r} \leq C \lambda _k^r \qquad (0 \leq r
\leq 1).
\end{equation}
According to~\eqref{eq.quasi},
$$e^{it (\partial_x \beta(x) \partial_x)} \phi_k = e^{it\lambda_k^2}
\phi_k + v$$ where $\|v\|_{L^\infty_{t,\text{loc}}( H^1(
\mathbb{R}))} = \mathcal{O} ( \lambda_k ^{- \infty})$.  Using the
Sobolev embedding $H^1 \rightarrow L^q$, we
can drop the contribution of $v$ in~\eqref{eq.nostrich}.  Using
H\"older inequality (and the fact that $\phi_k$ is supported in a ball
of radius $2^{-k}$), we obtain
$$ 1= \|\phi_k\|_{L^2}\leq C 2^{-k(q-2)/q}\|\phi_k\|_{L^q}$$ and
consequently, according to~\eqref{eq.sobo},
$$ \frac {\| e^{ it (\partial_x \beta(x) \partial_x)} \phi_k \|_{L^1(-
\varepsilon, \varepsilon); L^q}} { \|\phi_k \|_{H^r(
\mathbb{R})}}\geq\frac{ c 2^{k(q-2)/q}} {(k2^k)^r} \rightarrow +
\infty, \text{\ when\ } k\rightarrow + \infty.
$$ We now come back to the proof of Proposition~\ref{prop.quasi}. The
starting point is the interval instability of the Hill equation (see
for example~\cite{CoSp89}):
\begin{lemme} There exist $w, \alpha \in C^\infty$ such that
$$ w'' + \alpha w =0 \text{\ \ on\ \ } \mathbb{R},$$ $\alpha$ is 1-periodic
on $\mathbb{R}^+$ and $\mathbb{R}^-$, equal to $4\pi^2$ in a neighborhood
of $0$, and
$$| \alpha - 4 \pi^2| \leq 1,$$ $w(x)= p e^{-|x|}$ where $p$ is
1-periodic on $\mathbb{R}^+$ and $\mathbb{R}^-$, and $\| w \|_{L^2}=1$.
\end{lemme} 
Changing variables and setting
$$ y(x) = \int_0^x \alpha (s) ds, \qquad v(y) = w(x(y)), \qquad
\beta(y) = \alpha(x(y))$$ we get $$\frac{\partial} {\partial y} =
\alpha ^{-1} (x) \frac{\partial} {\partial x} \text{ and } (\partial_y
\beta(y) \partial_y +1) v=0$$ Denote by $$v^{\lambda, m}(y) =
v(\lambda (y-m)), \qquad \beta^{\lambda, m}(y) = \beta
(\lambda(y-m))$$ solutions of
\begin{align}
 (\partial_y \beta^{\lambda,m}(y) \partial_y +\lambda^2)
 v^{\lambda,m}&=0 \\
\label{eq.estmdec}
 |v^{\lambda, m} (y) | &\leq C e^{-\lambda |y-m|}\end{align}
 Consider
$\Psi_1 \in C^\infty_0 (]-1/4, 1/4[)$ equal to $1$ on $[-1/5,1/5]$,
$\Psi_2 \in C^\infty_0 (]-1/5, 1/5[)$ equal to $1$ on $[-1/6,1/6]$,
sequences $m_n= 2^{-n}$, $\lambda _n = n2^n$.

Using~\eqref{eq.estmdec}, we see that $v_n= v^{\lambda_n, m_n} (y)
\Psi_2 ( 2^n (y-m_n))$ is solution of
\begin{equation}
\label{eq.resc}
 (\partial_y \beta^{\lambda_n, m_n} (y) \partial_y + \lambda_n^2)v_n =
 \mathcal{O}(\lambda_ne^{-cn})_{H^1}
\end{equation}
Remark also that on the support of $v_n$, $\Psi_1( 2^n (y-m_n))=1$ and
consequently we can replace in~\eqref{eq.resc} $\beta^{\lambda_n,
m_n}(y)$ by $\beta_n(y)= \beta^{\lambda_n, m_n}(y)\Psi_1( 2^n
(y-m_n))$. Remark also that for $p\neq n$, the support of $v_n$ is
disjoint from the support of $\Psi_1( 2^p (y-m_p))$. Consequently, we
can replace in~\eqref{eq.resc} $\beta^{\lambda_n, m_n}(y)$ by
\begin{equation*}
\beta(y) = \sum_{n\in \mathbb{N}}\beta_n(y) + 4\pi (1-\sum_{n\in
\mathbb{N}} \Psi_1( 2^n (y-m_n))
\end{equation*}
(the last term being here only to ensure that $\beta(y) \geq 2\pi$).

To prove Proposition~\ref{prop.quasi}, it is now enough to show that
$\beta$ is in $ {W}^{s,1}$ for any $0\leq s <1$.  A direct calculation
shows that,
$$\| \beta_n\|_{{W}^{1,1}} \sim n, \qquad \| \beta_n\|_{L^1}\sim
2^{-n} \Rightarrow \| \beta_n\|_{\dot {W}^{s,1}}\leq C n^{s} 2^{-(1-s)n}$$
which implies that the series defining $\beta$ converges in $W^{s,1}$.
\bibliography{BV}

 \end{document}